\documentclass[twoside,11pt]{article}

\usepackage{graphicx}
\graphicspath{.}
\usepackage{amsmath,amssymb}

\RequirePackage{hyperref}
\RequirePackage{natbib}

\bibpunct{(}{)}{;}{a}{,}{,}

\hypersetup{
    colorlinks,%
    citecolor=black,%
    filecolor=black,%
    linkcolor=black,%
    urlcolor=black
}

\newcommand{\E}{\mathbb{E}}

\newcommand{\cB}{\mathcal{B}}
\newcommand{\cD}{\mathcal{D}}
\newcommand{\cM}{\mathcal{M}}
\newcommand{\cN}{\mathcal{N}}
\newcommand{\cP}{\mathcal{P}}
\newcommand{\cU}{\mathcal{U}}
\newcommand{\cX}{\mathcal{X}}
\newcommand{\un}{\mathbf{1}}
\newcommand{\dimX}{p}
\newcommand{\parenj}[1]{\mathopen{}\left( #1  \right) \mathclose{}}
\newcommand{\parenb}[1]{\bigl( #1  \bigr)}
\newcommand{\crochj}[1]{\mathopen{}\left[ #1 \right] \mathclose{}}

\newcommand{\sachantb}{\, \big| \,}
\newcommand{\Dpart}{\cD_{n_1}^1}
\newcommand{\Dlabels}{\cD_{n_2}^2}
\begin{document}

\title{Comments on: ``A Random Forest Guided Tour'' by G. Biau and E. Scornet}

\author{Sylvain Arlot         \and
        Robin Genuer 
}

\author{Sylvain Arlot  \\ \texttt{sylvain.arlot@math.u-psud.fr} \\
       Laboratoire de Math\'ematiques d'Orsay \\
Univ. Paris-Sud, CNRS, Universit\'e Paris-Saclay \\
91405 Orsay, France
       \and\ 
       Robin Genuer \\ \texttt{robin.genuer@isped.u-bordeaux2.fr} \\
Univ. Bordeaux, ISPED, Centre INSERM U-1219, \\
       INRIA Bordeaux Sud-Ouest, Equipe SISTM, \\
		146 rue L\'eo Saignat, \\
		F-33076 Bordeaux Cedex, 
		France}

\maketitle

\begin{abstract}
This paper is a comment on the survey paper by 
\citet{Bia_Sco:2016:TEST} about random forests. 
We focus on the problem of quantifying the impact of each ingredient  of random forests on their performance. 
We show that such a quantification is possible for a simple pure forest, 
leading to conclusions that could apply more generally. 
Then, we consider ``hold-out'' random forests, 
which are a good middle point between ``toy'' pure forests and Breiman's original random forests. 
\end{abstract}


We would like to thank G. Biau and E. Scornet for 
their clear and thought-provoking survey \citep{Bia_Sco:2016:TEST}. 
It shows that for understanding better random forests mechanisms, 
we must go beyond consistency results and 
\emph{quantify} the impact of each ingredient  of random forests 
on their performance. 

In this comment, we prove that such a quantification is possible for a simple pure forest, 
leading to conclusions that could apply more generally (Section~\ref{sec:agg-sub-ran}). 
Then, we consider in Section~\ref{sec:horf} ``hold-out'' random forests, 
which are a good middle point between ``toy'' pure forests and Breiman's original random forests.

\section{Aggregation, two kinds of randomization and resampling}
\label{sec:agg-sub-ran}

As clearly shown by \citet{Bia_Sco:2016:TEST}, 
two key ingredients of random forests are the \emph{aggregation} process 
(several trees are combined to get the final estimator) 
and the \emph{diversity} of the trees that are aggregated. 
We can distinguish two sources of diversity: 
\begin{itemize}
\item randomization of the \emph{partition} $\cP_{\mathrm{final}}$, 
\item randomization of the \emph{labels}, that is, of the predicted value in each cell 
of $\cP_{\mathrm{final}}$, given $\cP_{\mathrm{final}}$. 
\end{itemize}
In Breiman's random forests, resampling acts on the two kinds of randomization, 
while the choice of $\cM_{\mathrm{try}}$ at each node of each tree acts on the randomization of the partitions only. 
In purely random forests, partitions are built independently from the data, 
so resampling (if any) only acts on the randomization of the labels. 
Therefore, the role of each kind of randomization is easier to understand, 
and to quantify separately, 
for purely random forests. 
In this section, we propose to do so for the one-dimensional toy forest 
introduced by \citet[Section~4]{Arl_Gen:2014}. 

\subsection{Toy forest}
\label{sec:toy}

We first define the toy forest model, 
assuming that $\cX=[0,1]$. 
Let $2 \leq k \leq a \leq n$ be some integers. 
The partition associated to each tree is given by 
\[ 
\left[ 0, \frac{1-T}{k} \right) , \, 
\left[ \frac{1-T}{k} , \frac{2-T}{k}  \right) , \, 
\ldots, 
\left[ \frac{k-1-T}{k} , \frac{k-T}{k}  \right) , \, 
\left[ \frac{k-T}{k} , 1  \right) 
\]
with $T \sim \cU([0,1])$. 
Then, for each tree, a subsample $(X_i,Y_i)_{i \in I}$ 
of size $a$ is chosen (uniformly over the subsamples of size $a$), 
independently from $T$. 
The tree estimator is defined as usual: 
the predicted value at $\mathbf{x}$ 
is the average of the $Y_i$ such that 
$X_i$ belongs to $A_n(x;T)$, 
the cell of the partition that contains $\mathbf{x}$. 
Finally, the forest estimator is obtained by 
averaging $M \geq 1$ trees, 
where the trees are independent conditionally to $\cD_n$. 

\medbreak 

Assume that $X_i \sim \cU([0,1])$, 
the $(X_i,Y_i)$ are independent with the same distribution, 
the noise-level $\E [ (Y-m(X))^2 \, | \, X ]=\sigma^2$ is constant, 
$m$ is of class $\mathcal{C}^3$, 
$n \gg 1$ and $a \gg k \log(n)$. 
Let $\mathbf{x} \in [k^{-1}, 1-k^{-1}]$ (to avoid
border effects, see \citet{Arl_Gen:2014} for details) be fixed.
Then, Table~\ref{tab.toy} provides the order of magnitude of 
\[ 
\E\Bigl[ \bigl( m_{M,n}^{\mathrm{toy}} (\mathbf{x}) - m(\mathbf{x}) \bigr)^2  \Bigr]
\enspace , 
\]
the quadratic risk at $\mathbf{x}$ of a toy forest 
$m_{M,n}^{\mathrm{toy}}$ with $M$ trees, 
in various situations: 
with or without aggregation ($M=+\infty$ or $M=1$), 
with or without randomization of the partitions (we remove randomization of partitions by putting $T=0$ for all trees). 
Subsampling can also be removed by taking $a=n$ in all formulas. 
Note that in Table~\ref{tab.toy}, 
the risk is written as the sum of an approximation error 
and an estimation error, see Appendix~\ref{sec:approx-estim}. 
The main lines of the proof are given in Appendix~\ref{sec:toy-proof}. 
\begin{table} 
\caption{Order of magnitude of the quadratic risk at $\mathbf{x}$ (approximation error + estimation error) 
for the one-dimensional toy random forest 
combined with subsampling without replacement; 
$a$ denotes the subsample size, 
$1/k$ is the size of each cell in a partition, 
see Appendix~\ref{sec:toy-proof} for precise assumptions and approximations. 
The impact of subsampling can be quantified by comparing $a=n$ with $a<n$ in each formula. 
\label{tab.toy}}
\begin{center}
\begin{tabular}{ccc}
\hline\noalign{\smallskip}
& Single tree & Infinite forest 
\\
\noalign{\smallskip}\hline\noalign{\smallskip}
%
\begin{minipage}{2cm}\centerline{No randomization} 
\centerline{of partitions} \end{minipage} 
& $\displaystyle \frac{c_1(m,\mathbf{x})}{k^{2}} + \frac{\sigma^{2\vphantom{^2}} k}{a}$ & $\displaystyle\frac{c_1(m,\mathbf{x})}{k^2\vphantom{_3}} + \frac{\sigma^2 k}{n}$
\\
\noalign{\smallskip}\hline\noalign{\smallskip}
\begin{minipage}{2cm} \centerline{Randomization} 
\centerline{of partitions}  \end{minipage} & $\displaystyle\frac{c_1(m,\mathbf{x})}{k^2} + \frac{\sigma^{2\vphantom{^2}} k}{a}$ & $\displaystyle \frac{c_2(m,\mathbf{x})}{k^4\vphantom{_3}} + \frac{2 \sigma^2 k}{3 n}$
\\
\noalign{\smallskip}\hline
\end{tabular} 
\[ 
\text{where} \qquad 
c_1(m,\mathbf{\mathbf{x}}) = \frac{m'(\mathbf{x})^2}{12}
\qquad \text{and} \qquad
c_2(m,\mathbf{\mathbf{x}}) = \frac{m^{\prime\prime}(\mathbf{x})^2}{144}
\enspace . 
\]
\end{center}
\end{table}

\medbreak

%
Table~\ref{tab.toy} allows to quantify the impact 
of aggregation, randomization of the partitions, 
randomization of the labels---and their combinations---on the performance of toy forests: 
\begin{itemize}
\item \emph{aggregation}: 
comparing the two columns of Table~\ref{tab.toy} 
shows that aggregation always improve the performance (which is true for any forest model, by Jensen's inequality). 
The improvement can be huge: 
when partitions and labels are randomized, 
$a \ll n$ and $k \gg 1$,  
both approximation and estimation errors decrease by an order of magnitude. 

\item \emph{randomization of the partitions}: 
comparing the two lines of Table~\ref{tab.toy} 
shows that randomizing the partitions strongly improves the performance of the infinite forest 
(there is no change for a single tree, of course). 
The approximation error decreases by an order of magnitude, 
as previously showed by \citet{Arl_Gen:2014}. 
The estimation error also decreases---as showed by \citet{Gen:2012} for another pure forest---, 
but by a factor $3/2$ only. 

\item \emph{randomization of the labels}: 
comparing $a=n$ with $a \ll n$ in the formulas of Table~\ref{tab.toy} 
shows the influence of subsampling, 
which is the only randomization mechanism for the labels. 
Single trees perform worse when $a < n$ 
(as expected since the sample size is lowered). 
The performance of infinite forests does not change with 
subsampling, 
which might seem a bit surprising given several results 
mentioned by \citet{Bia_Sco:2016:TEST}. 
This phenomenon corresponds to the fact that subagging does not improve a stable estimator \citep{Buh_Yu:2002}, 
and that a regular histogram is stable.
Section~\ref{sec:agg-sub-ran:discuss} below explains why there is no contradiction with the random forests literature.  
\end{itemize}

\subsection{Discussion}
\label{sec:agg-sub-ran:discuss}
%
%
Section~\ref{sec:toy} sheds some light on 
previous theoretical results on random forests, 
and suggests a few conjectures which deserve to be investigated. 

\paragraph{Parametrization of the trees.} 
%
%
The end of Section~3.1 of Biau and Scornet's survey  
might seem contradictory with the above results for the toy forest. 
According to most papers in the literature, 
``random forests reduce the estimation error of a single tree, 
while maintaining the same approximation error''. 
Moreover, an infinite forest can be consistent even when a single tree (grown with a sample of size $n$) is not. 
Table~\ref{tab.toy} precisely shows the opposite situation: 
the estimation error is almost the same for a single tree and for an infinite forest, 
while the approximation error is dramatically smaller for an infinite forest. 
In addition, when an infinite forest is consistent, $1 \ll k \ll n$ 
hence a single tree trained with $a=n$ points is also consistent. 

%
The point is that these results and ours 
consider different parametrizations of the trees. 
In Section~\ref{sec:toy}, trees are parametrized by the number of leaves $k+1$; 
so, when comparing a tree with a forest, 
we think fair to compare 
(i) a tree of $k+1$ leaves trained with $n$ data points, 
with (ii) a forest where each tree has $k+1$ leaves and is trained with $a$ data points. 
In the literature, trees are often parametrized by 
the number $\mathtt{nodesize}$ of data points per cell. 
Then, comparisons are done between 
(i) a tree of $\approx n/\mathtt{nodesize}$ leaves trained with $n$ data points, 
and (ii) a forest where each tree has $\approx a/\mathtt{nodesize}$ leaves and 
is trained with $a$ data points. 
So, if we take $k \approx a/\mathtt{nodesize}$, the two approaches consider 
(approximately) the same forest (ii), 
but the reference trees (i) are quite different.

\medbreak 

%
We do not mean that one of these two parametrizations 
is definitely better than the other: 
$\mathtt{nodesize}$ is a natural parameter for Breiman's random forests, 
while toy forests are naturally parametrized by their number of leaves. 
Nevertheless, one must keep in mind that any comparison between a forest and a single tree trained with the full sample 
\emph{does depend on the parametrization}. 

The parametrization by $\mathtt{nodesize}$ 
can also hide some difficulties. 
For the toy forest model, $k$ has to be chosen, 
and this is not an easy task. 
One could think that this problem is solved by taking the $\mathtt{nodesize}$ 
parametrization with, say, $\mathtt{nodesize}=1$. 
This is wrong because we then have to choose the subsample size $a$, 
which is equivalent to the original problem since $k \approx a/\mathtt{nodesize}$.

\paragraph{What about Breiman's forests?} 
Section~\ref{sec:toy} suggests that for the toy forest, 
the most important ingredient in the tree diversity 
is the randomization of the partitions. 
We conjecture that this holds true for general random forests. 

Nevertheless, we do not mean that the resampling step should 
always be discarded, 
since for Breiman's random forests (for instance), 
resampling also acts on the randomization of the partitions. 
A key open problem is to quantify the relative roles of 
resampling and of the choice of $\cM_{\mathrm{try}}$ 
on the randomization of the partitions. 
Section~\ref{sec:horf} below shows that  
``hold-out random forests'' can be a good playground 
for such investigations. 

\paragraph{Bootstrap or subsampling?} 
Another important question is the choice between 
the bootstrap and subsampling, 
which remains an open problem according to 
\citet{Bia_Sco:2016:TEST}. 

We conjecture that Table~\ref{tab.toy} is also valid for the $a$ out of $n$ bootstrap, 
which would mean that bootstrap and subsampling are fully 
equivalent with respect to the randomization of the \emph{labels}, 
with no impact on the performance of pure forests. 
Assuming this holds true, subsampling (with or without replacement) with $a \ll n$ remains 
interesting for reducing the computational cost---Table~\ref{tab.toy} only requires 
that $a \gg k \log(n)$. 

A key open problem remains: 
compare bootstrap and subsampling with respect to the 
randomization of the \emph{partitions} 
for Breiman's random forests. 
Again, ``hold-out random forests'' described 
in Section~\ref{sec:horf} should be a good starting point 
for such a comparison.

\section{Hold-out random forests}
\label{sec:horf}

We now consider a more complex forest model, called hold-out random forests, which is close to Breiman's random forests while being simpler to analyze. 
%
%
Hold-out random forests have been proposed by \citet[Section~3]{Bia:2012} 
and appear in the experiments of \citet[Section~7]{Arl_Gen:2014}. 

\subsection{Definition}

%
Hold-out random forests can be defined as follows. 
First, the data set $\cD_n$ is split, once and for all, 
into two subsamples $\Dpart$ and $\Dlabels$, 
of respective sizes $n_1$ and $n_2$, 
satisfying $n = n_1+n_2$. 
This split is done independently from $\cD_n$. 
Then, conditionally to $(\Dpart,\Dlabels)$, 
the $M$ trees are built independently as follows. 
The partition associated with the $j$-th tree is built as 
for Breiman's random forests with $\Dpart$ as data set. 
In other words, it is the partition $\cP_{\mathrm{final}}$ 
defined by \citet[Algorithm~1]{Bia_Sco:2016:TEST}  
with training set $\Dpart$ as an input. 
The $j$-th tree estimate at  $\mathbf{x}$ is defined by 
\[
m^{\mathrm{hoRF}}_n \parenb{ \mathbf{x} ; \Theta_j , \Dpart , \Dlabels }
:= \sum_{(X_i,Y_i) \in \Dlabels} \frac{\un_{X_i \in A_{n_1}(\mathbf{x} ; \Theta_j , \Dpart) } Y_i}{N_{n_2}( \mathbf{x} ; \Theta_j, \Dpart, \Dlabels) }
\]
where $A_{n_1}(\mathbf{x} ; \Theta_j , \Dpart)$ 
is the cell of this partition that contains $\mathbf{x}$ 
and $N_{n_2}( \mathbf{x} ; \Theta_j, \Dpart, \Dlabels)$ 
is the number of points $(X_i,Y_i) \in \Dlabels$ such that $X_i \in A_{n_1}(\mathbf{x} ; \Theta_j , \Dpart)$. 
Finally, the hold-out forest estimate is defined by 
\[ 
m^{\mathrm{hoRF}}_{M,n} \parenb{\mathbf{x} ; \Theta_{1 \ldots M} , \Dpart , \Dlabels }
= \frac{1}{M} \sum_{j=1}^M m^{\mathrm{hoRF}}_n \parenb{ \mathbf{x} ; \Theta_j , \Dpart , \Dlabels }
\, . 
\]

\medbreak

%
In the definition of $m^{\mathrm{hoRF}}_{M,n}$, 
the building of the partitions depends on the same parameters as Breiman's random forests: 
$\mathtt{mtry}$, $\mathtt{nodesize}$, 
the fact that resampling can be done with or without replacement, 
and the resample size $a_{n_1}$. 
It is also possible to add another resampling step when assigning labels to the leaves of each tree with $\Dlabels$; 
we do not consider it here since Section~\ref{sec:agg-sub-ran} suggests that this would not change much the performance (at least for forests). 

\medbreak

%
A key property of hold-out random forests is that they 
are \emph{purely random}: 
the subsamples $\Dpart$ and $\Dlabels$ are independent, 
hence the partition associated with the trees are independent from $\Dlabels$. 
Note however that each partition still depends on some 
data (through $\Dpart$), 
hence it could adapt to some features of the data, 
such as the ``sparsity'' of the regression function, 
the non-uniformity over $\cX$ of the smoothness of $m$,  
or the non-uniformity of the distribution of $X$. 
Therefore, hold-out random forests can capture much 
of the complexity of Breiman's random forests, 
while being easier to analyze since they are purely random. 
%
%
In particular, we can apply the results proved in Appendix~\ref{sec:approx-estim} (where our $\Dlabels$ is written $\cD_n$, and our $\Dpart$ is hidden in the relationship between $\Theta_j$ and the partition of the $j$-th tree). 
Then, the quadratic risk of $m^{\mathrm{hoRF}}_{M,n}$ 
can be (approximately) decomposed into the sum of an 
approximation error and an estimation error, 
and these two terms can be studied separately. 
For instance, the results of \citet{Arl_Gen:2014} can 
be applied in order to analyze the approximation error. 

%
For now, we only study the behaviour of these two 
terms in a short numerical experiment. 
The results are summarized by Table~\ref{tab.hoRF}, 
where estimated values of the approximation and estimation errors 
are reported as a function of $n_2$, $\sigma^2$ and of the parameters 
of the partition building process ($k$, $\mathtt{mtry}$ and bootstrap). 
\begin{table} 
\caption{Numerical estimation of the quadratic risk (approximation error + estimation error) 
for the hold-out random forest; 
$k$ is the number of cells in the partition; 
$\cX = [0,1]^p$ with $p=5$. 
There is no randomization of the labels.  
\label{tab.hoRF}}
\begin{center}
\begin{tabular}{ccc}
\hline\noalign{\smallskip}
& Single tree & Large forest 
\\
\noalign{\smallskip}\hline\noalign{\smallskip}
%
\begin{minipage}{2cm}\centerline{No bootstrap} 
\centerline{$\mathtt{mtry}=p$} \end{minipage} 
& $\displaystyle \frac{0.13}{k^{0.17}} + \frac{1.04 \sigma^{2\vphantom{^2}} k}{n_2}$
& $\displaystyle \frac{0.13}{k^{0.17}} + \frac{1.04 \sigma^{2\vphantom{^2}} k}{n_2}$
\\
\noalign{\smallskip}\hline\noalign{\smallskip}
\begin{minipage}{2cm} \centerline{Bootstrap}  
\centerline{$\mathtt{mtry}=p$} \end{minipage}
& $\displaystyle \frac{0.14}{k^{0.17}} + \frac{1.06 \sigma^{2\vphantom{^2}} k}{n_2}$
& $\displaystyle \frac{0.15}{k^{0.29}} + \frac{0.08 \sigma^{2\vphantom{^2}} k}{n_2}$
\\
\noalign{\smallskip}\hline
\begin{minipage}{2cm} \centerline{No bootstrap}  
\centerline{$\mathtt{mtry}=\lfloor p/3 \rfloor$} \end{minipage}
& $\displaystyle \frac{0.23}{k^{0.19}} + \frac{1.01 \sigma^{2\vphantom{^2}} k}{n_2}$
& $\displaystyle \frac{0.06}{k^{0.31}} + \frac{0.06 \sigma^{2\vphantom{^2}} k}{n_2}$
\\
\noalign{\smallskip}\hline
\begin{minipage}{2cm} \centerline{Bootstrap}  
\centerline{$\mathtt{mtry}=\lfloor p/3 \rfloor$} 
\end{minipage}
& $\displaystyle \frac{0.25}{k^{0.20}} + \frac{1.02 \sigma^{2\vphantom{^2}} k}{n_2}$
& $\displaystyle \frac{0.06}{k^{0.34}} + \frac{0.05 \sigma^{2\vphantom{^2}} k}{n_2}$
\\
\noalign{\smallskip}\hline
\end{tabular} 
\end{center}
\end{table}
Detailed information about this experiment can 
be found in Appendix~\ref{sec:app-expe}. 
Let us emphasize here that we consider a single data generation setting, 
hence these results must be interpreted with care.

\subsection{Discussion}
Based on Table~\ref{tab.hoRF} and our experience about 
Breiman's random forests, we can make the following comments. 

\paragraph{Choice of $\mathtt{mtry}$.}

As illustrated in Table~\ref{tab.hoRF}, choosing 
$\mathtt{mtry} = \lfloor p/3 \rfloor$ 
instead of $p$ 
decreases the risk of an infinite forest.
When there is no bootstrap, the performance gain is
significant and the reason is that it is the only source of randomization of
partitions. 
But, even in presence of bootstrap, it allows to slightly reduce
the approximation error. 
In the same experiments with $p=10$, 
the gain of decreasing $\mathtt{mtry}$
in the bootstrap case is larger 
(see the supplementary material).

Our belief is that when there is some bootstrap, 
the additional randomization given by taking $\mathtt{mtry} < p$
can reduce the risk in some cases, where typically $n\geq p$
(which holds true in our experiments).
This is supported by the experiments of \citet[Section~2]{Gen_Pog_Tul:2008},  
where small values of \texttt{mtry} give significantly lower risk than
$\mathtt{mtry}=p$ for some classification problems. 
For regression, \citet[Section~2]{Gen_Pog_Tul:2008} obtain 
similar performance when decreasing \texttt{mtry}, 
which is consistent with Table~\ref{tab.hoRF} since 
these experiments are done in the bootstrap case.

When $n \ll p$ and only a small proportion of the coordinates 
of $\mathbf{x}$ are informative, 
we conjecture that the optimal $\mathtt{mtry}$ is close to $p$ 
(provided that there is some bootstrap step for randomizing the partitions). 
Indeed, if \texttt{mtry} is significantly smaller than $p$, 
then, the probability to choose at least one informative coordinate
in $\cM_{\mathrm{try}}$ is not close to~$1$, 
hence the randomization of the partitions might be too strong. 

\paragraph{Bootstrap, $\mathtt{mtry}$ and randomization of the partitions.}
When $\mathtt{mtry}=p$, according to Table~\ref{tab.hoRF}, 
the bootstrap helps to significantly reduce
the risk, compared with the no randomization case. 
Overall, we get a 
significantly smaller risk when there is \emph{at least one}  
source of randomization of the partitions. 

Comparing the three combinations of parameters 
(bootstrap, $\mathtt{mtry} <p$, or both) for which 
the partitions are randomized is more difficult: 
the differences observed in Table~\ref{tab.hoRF} might not 
be significant. 
Nevertheless, Table~\ref{tab.hoRF} suggests that the lowest 
risk might be obtained when two sources of randomization 
are present ($\mathtt{mtry} <p$ and bootstrap). 
And if we have to choose only one source of 
randomization, it seems that randomizing with $\mathtt{mtry} 
= \lfloor p/3 \rfloor$ only yields a smaller risk than 
bootstrapping only.

\appendix 
\section{Approximation and estimation errors}
\label{sec:approx-estim}

We state a general decomposition of the risk of a forest 
having the $X$-property (that is, when partitions are built independently from $(Y_i)_{1 \leq i \leq n}$), 
that we need for proving the results of Section~\ref{sec:agg-sub-ran}, 
but can be useful more generally. 
We assume that $\E[Y_i^2]<+\infty$ for all $i$. 

%
For any random forest $m_{M,n}$ having the $X$-property, 
following \citet[Sections~2 and~3.2]{Bia_Sco:2016:TEST}, 
we can write 
\begin{gather} 
\notag 
m_{M,n} (\mathbf{x};\Theta_{1 \ldots M}, \cD_n) 
= \sum_{i=1}^n W_{ni}(\mathbf{x}) Y_i
\\
\label{eq:Wni} 
\text{where} \qquad 
W_{ni}(\mathbf{x}) 
= W_{ni}(\mathbf{x} ;\Theta_{1 \ldots M} , X_{1 \ldots n}) 
= \frac{1}{M} \sum_{j=1}^M \frac{C_i(\Theta_j) \un_{X_i \in A_n(\mathbf{x}; \Theta_j ; X_{1 \ldots n})}}{N_n(\mathbf{x}; \Theta_j ; X_{1 \ldots n})}
\enspace , 
\end{gather}
$C_i(\Theta_j)$ is the number of times $(X_i,Y_i)$ appears in the $j$-th resample, 
$A_n(\mathbf{x}; \Theta_j ; X_{1 \ldots n})$ is the cell containing 
$\mathbf{x}$ in the $j$-th tree, 
and 
\[ 
N_n(\mathbf{x}; \Theta_j ; X_{1 \ldots n})
= 
\sum_{i=1}^n C_i(\Theta_j) \un_{X_i \in A_n(\mathbf{x}; \Theta_j ; X_{1 \ldots n})}
\enspace . 
\]
%
%
Now, let us define 
\begin{align*} 
m^{\star}_{M,n} (\mathbf{x};\Theta_{1 \ldots M}, X_{1 \ldots n}) 
&= \E \Bigl[ m_{M,n} (\mathbf{x};\Theta_{1 \ldots M}, \mathcal{D}_n)  \sachantb X_{1 \ldots n}, \Theta_{1 \ldots M} \Bigr]
\\
&= \sum_{i=1}^n W_{ni}(\mathbf{x};\Theta_{1 \ldots M} , X_{1 \ldots n}) m(X_i)
\\
\text{and} \qquad  
%
\overline{m}^{\star}_{M,n} (\mathbf{x};\Theta_{1 \ldots M}) 
&= 
\E \Bigl[ m^{\star}_{M,n} (\mathbf{x};\Theta_{1 \ldots M}, X_{1 \ldots n})  \sachantb \Theta_{1 \ldots M} \Bigr]
\enspace . 
\end{align*}
By definition of the conditional expectation, 
we can decompose the risk of $m_{M,n}$ at $\mathbf{x}$ into three terms 
\begin{equation}
\label{eq:decomp}
\begin{split}
&\hspace*{-1cm}\E\Bigl[ \bigl( m_{M,n}(\mathbf{x}) - m (\mathbf{x}) \bigr)^2 \Bigr]
= 
\underbrace{
\E\Bigl[ \bigl( \overline{m}^{\star}_{M,n} (\mathbf{x}) - m(\mathbf{x}) \bigr)^2 \Bigr]
}_{A = \text{approximation error}}
\\  \hspace*{1cm}
&+ 
\underbrace{
\E\Bigl[ \bigl( m^{\star}_{M,n} (\mathbf{x}) - \overline{m}^{\star}_{M,n}(\mathbf{x}) \bigr)^2 \Bigr]
}_{\Delta}
+ 
\underbrace{
\E\Bigl[ \bigl( m_{M,n}(\mathbf{x}) - m^{\star}_{M,n} (\mathbf{x}) \bigr)^2 \Bigr]
}_{E = \text{estimation error}}
\enspace  . 
\end{split}
\end{equation}
In the fixed-design regression setting (where the $X_i$ are deterministic), 
$A$ is called approximation error, $\Delta=0$, and $E$ is called estimation error. 
Things are a bit more complicated in the random-design setting---when $(X_i,Y_i)_{1 \leq i \leq n}$ are independent and identically distributed---since $\Delta  \neq 0$ in general. 
Up to minor differences related to how $m_{n}$ is defined on empty cells, 
$A$ is still the approximation error, 
and the estimation error is $\Delta + E$. 

Let us finally assume that 
$(X_i,Y_i)_{1 \leq i \leq n}$ are independent 
and define 
\[ 
\sigma^2(X_i) = 
\E\crochj{ \bigl( m(X_i) - Y_i \bigr)^2 \sachantb X_i }
\enspace . \]
Then, since the weights $W_{ni}(\mathbf{x})$ only depend on $\cD_n$ through $X_{1 \ldots n}$, 
we have the following formula for the estimation error
\begin{align}
\notag 
E = \E\crochj{ \mathopen{} \left( \sum_{i=1}^n W_{ni}(\mathbf{x}) \bigl( m(X_i) - Y_i \bigr) \right)^2 \mathclose{} } 
&= 
\E\crochj{ \sum_{i=1}^n W_{ni}(\mathbf{x})^2 \sigma^2(X_i) } 
\enspace . 
\end{align}
For instance, in the homoscedastic case, 
$\sigma^2(X_i) \equiv \sigma^2$ and  
\begin{align}
\label{eq:est-err}
E = \E\crochj{ \mathopen{} \left( \sum_{i=1}^n W_{ni}(\mathbf{x}) \bigl( m(X_i) - Y_i \bigr) \right)^2 \mathclose{} } 
&= 
\sigma^2 \E\crochj{ \sum_{i=1}^n W_{ni}(\mathbf{x})^2 } 
\enspace . 
\end{align}

\section{Analysis of the toy forest: proofs}
\label{sec:toy-proof}

We prove the results stated in Section~\ref{sec:agg-sub-ran} for the one-dimensional 
toy forest. 

Since the toy forest is purely random, all results of Appendix~\ref{sec:approx-estim} apply, 
with $\Theta = (T,I)$ and $C_i(\Theta) = \un_{i \in I}$. 
It remains to compute the three terms of Eq.~\eqref{eq:decomp}. 

%
Since we assume $m$ is of class $\mathcal{C}^3$, 
we can use the results of \citet[Section~4]{Arl_Gen:2014} 
for the approximation error $A$ 
(up to minor differences in the definition of $\overline{m}^{\star}_{M,n} (\mathbf{x})$, due to event where $A_n(\mathbf{x};\Theta)$ is empty, which has a small probability since $a \gg k$). 
We assume that $m'(\mathbf{x}) \neq 0$ and $m''(\mathbf{x}) \neq 0$ 
for simplicity, 
so the quantities appearing in Table~\ref{tab.toy} 
indeed provide the order of magnitude of~$A$. 

%
The middle term $\Delta$ in decomposition \eqref{eq:decomp} 
is negligible in front of $E$ for a single tree, 
which can be proved using results from \citet{Arl:2008a}, 
as soon as $m'(\mathbf{x}) / k \ll \sigma$ and $a \gg k$. 
We assume that it can also be neglected for an infinite forest. 

%
For the estimation error, we can use Eq.~\eqref{eq:est-err} 
and the following arguments. 
First, for every $i \in \{1, \ldots, n\}$, $X_i$ belongs to $A_n(\mathbf{x};\Theta)$ 
with probability $1/k$. 
Combined with the subsampling process, 
we get that 
\[
N_n(\mathbf{x};\Theta; X_{1 \ldots n}) \sim \cB\parenj{ n , \frac{a}{n k} }
\]
is close to its expectation $a/k$ with probability almost one 
if $a/k \gg \log(n)$. 
Assuming that this holds simultaneously for a huge fraction of the subsamples, 
we get the approximation 
\begin{align}
\notag 
W_{ni}^{\mathrm{toy}}(\mathbf{x}) 
&= 
\frac{1}{M} \sum_{j=1}^M \frac{\un_{i \in I_j} \un_{X_i \in A_n(\mathbf{x}; \Theta_j)}}{N_n(\mathbf{x}; \Theta_j ; X_{1 \ldots n})}
\\
&\approx 
\frac{k}{a}
\frac{1}{M} \sum_{j=1}^M \un_{i \in I_j} \un_{X_i \in A_n(\mathbf{x}; \Theta_j)}
=:
\widetilde{W}_{ni}^{\mathrm{toy}}(\mathbf{x}) 
\enspace . 
\label{eq:toy:approx-W}
\end{align}

Now, we note that conditionally to $X_{1 \ldots n}$, 
the variables $\un_{i \in I_j} \un_{X_i \in A_n(\mathbf{x}; \Theta_j)}$, 
$j=1, \ldots, M$ are independent and follow a Bernoulli distribution with 
the same parameter 
\[
\frac{a}{n} \times \bigl( 1-k|X_i-x| \bigr)_+ 
\enspace . 
\]
Therefore, 
\begin{align*}
\E\crochj{ \widetilde{W}_{ni}^{\mathrm{toy} }(\mathbf{x})^2 \, \big| \, X_{1 \ldots n}} 
&= 
\frac{k^2 }{ n a} 
\crochj{ 
\parenj{ 1 - \frac{1}{M} } \frac{a}{n}
\Bigl( \bigl( 1-k|X_i-x| \bigr)_+ \Bigr)^2
+ 
\frac{1}{M}  \bigl( 1-k|X_i-x| \bigr)_+
}
\\
\text{hence} \quad 
\E\crochj{ \widetilde{W}_{ni}^{\mathrm{toy} }(\mathbf{x})^2 } 
&= 
\frac{k}{n a} 
\crochj{ 
\parenj{ 1 - \frac{1}{M} } \frac{2 a}{3 n}
+ 
\frac{1}{M}
}
\enspace . 
\end{align*}
By Eq.~\eqref{eq:est-err}, this ends the proof 
of the results in the bottom line of Table~\ref{tab.toy}. 

Similar arguments apply for justifying the top line of Table~\ref{tab.toy}, 
where $T_j=0$ almost surely. 

\medbreak

Note that we have not given a full rigorous proof of the results shown in Table~\ref{tab.toy}, 
because of the approximation~\eqref{eq:toy:approx-W} and
of the term $\Delta$ that we have neglected. 
We are convinced that the parts of the proof 
that we have skipped might only require to add some technical assumptions, 
which would not help to reach our goal of understanding better random forests in general.

\section{Details about the experiments}
\label{sec:app-expe}

This section describes the experiments whose results are shown 
in Section~\ref{sec:horf}. 

\paragraph{Data generation process.}
We take $\cX = [0,1]^{\dimX}$, with $p \in \{ 5, 10 \}$. 
Table~\ref{tab.hoRF} only shows the results for $p=5$. 
Results for $p=10$ are shown in supplementary material. 

The data $(X_i,Y_i)_{1 \leq i \leq n_1 + n_2}$ are independent with 
the same distribution: 
$X_i \sim \cU([0,1]^{\dimX})$, 
$Y_i = m(X_i) + \varepsilon_i$ with $\epsilon_i \sim \cN(0,\sigma^2)$ 
independent from $X_i$, $\sigma^2 = 1/16$, 
and the regression function $m$ is defined by 
  \[ m : \mathbf{x} \in [0,1]^{\dimX} \mapsto \mathbf{1/10} \times \left[ 10 \sin(\pi x_1 x_2) + 20 (x_3 -
  0.5)^2 + 10 x_4 + 5 x_5 \right] 
  \, . \]
The function $m$ is proportional to the \textbf{Friedman1} function 
which was introduced by \citet{Fri:1991}.
Note that when $\dimX>5$, $m$ only depends on the $5$ first
coordinates of $\mathbf{x}$.

Then, the two subsamples are defined by 
$\Dpart = (X_i,Y_i)_{1 \leq i \leq n_1}$ and $\Dlabels = (X_i,Y_i)_{n_1 + 1 \leq i \leq n_1+n_2}$. 

We always take $n_1 = 1\,280$ and $n_2 = 25\,600$. 

\paragraph{Trees and forests.}
For each $k \in \{2^5, 2^6, 2^7, 2^8\}$, 
each experimental condition (bootstrap or not, $\mathtt{mtry}=p$ or $\lfloor p/3 \rfloor$), 
we build some hold-out random trees and forests 
as defined in Section~\ref{sec:horf}. 
These are built with the \texttt{randomForest} R package
\citep{Lia_Wie:2002,R_Core:2014}, with appropriate parameters
($k$ is controlled by \texttt{maxnodes}, while $\mathtt{nodesize}=1$). 

Resampling within $\Dpart$ (when there is some resampling) 
is done with a bootstrap sample of size $n_1$ 
(that is, with replacement and $a_{n_1} = n_1$). 

``Large'' forests are made of $M=k$ trees, a number of trees suggested by \citet{Arl_Gen:2014}.

\paragraph{Estimates of approximation and estimation error.}
Estimating approximation and estimation errors (as defined by Eq.~\eqref{eq:decomp}) 
requires to estimate some expectations over $\Theta$ 
(which includes the randomness of $\Dpart$ as well as the randomness of the 
choice of bootstrap subsamples of $\Dpart$ and of the repeated choices 
of a subset $\cM_{\mathrm{try}}$). 
This is done with a Monte-Carlo approximation, 
with $500$ replicates for trees 
and $10$ replicates for forests. 
This number might seem small, but we observe that large forests are quite 
stable, hence expectations can be evaluated precisely from a small number of replicates.

We estimate the approximation error (integrated over $\mathbf{x}$) as follows. 
For each partition that we build, we compute the corresponding 
``ideal'' tree, which maps each piece of the partition 
to the average of $m$ over it 
(this average can be computed almost exactly from the definition of $m$). 
Then, to each forest we associate the ``ideal'' forest 
$\overline{m}^{\star}_{M,n}$ which is the average of the ideal trees. 
We can thus compute 
$( \overline{m}^{\star}_{M,n} (\mathbf{x}) - m(\mathbf{x}) )^2$ 
for any $\mathbf{x} \in \cX$, 
and estimate its expectation with respect to $\Theta$. 
Averaging these estimates over $1000$ uniform random points $\mathbf{x} \in \cX$ 
provides our estimate of the approximation error. 

We estimate the estimation error (integrated over $\mathbf{x}$) from Eq.~\eqref{eq:est-err}; 
since $\sigma^2$ is known, we focus on the remaining term. 
Given some hold-out random forest, 
for any $\mathbf{x} \in \cX$ and $i \in \{1, \ldots, n\}$, 
we can compute
\begin{equation} 
\notag
W_{ni}(\mathbf{x}) 
= \frac{1}{M} \sum_{j=1}^M \sum_{(X_i,Y_i) \in \Dlabels} \frac{\un_{X_i \in A_{n_1}(\mathbf{x} ; \Theta_j , \Dpart) }}{N_{n_2}( \mathbf{x} ; \Theta_j, \Dpart, \Dlabels) }
\enspace .
\end{equation}
Then, averaging $\sum_i W_{ni}(\mathbf{x})^2$ over 
several replicate trees/forests 
and over $1\,000$ uniform random points $\mathbf{x} \in \cX$, 
we get an estimate of the estimation error (divided by $\sigma^2$).

\paragraph{Summarizing the results in Table~\ref{tab.hoRF}.}
Given the estimates of the (integrated) approximation and estimation errors 
that we obtain for every $k \in \{2^5, 2^6, 2^7, 2^8\}$, 
we plot each kind of error as a function of $k$ 
(in $\mathrm{log}_2$-$\mathrm{log}_2$ scale for the approximation error), 
and we fit a simple linear model (with an intercept). 
The estimated parameters of the model directly give the results shown in 
Table~\ref{tab.hoRF} 
(in which the value of the intercept for the estimation error is omitted for simplicity). 
The corresponding graphs are shown in supplementary material.

\section*{Acknowledgements}
The research of the authors is partly supported by 
the French Agence Nationale de la Recherche (ANR 2011 BS01 010 01 projet Calibration). 
S. Arlot is also partly supported by Institut des Hautes \'Etudes Scientifiques (IHES, Le Bois-Marie, 35, route de Chartres,
91440 Bures-Sur-Yvette, France).

\bibliography{biblio}

\clearpage

\section{Supplementary Material}
\label{sec:sup-mat}

\begin{table} 
\caption{Numerical estimation of the quadratic risk (approximation error + estimation error) 
for the hold-out random forest; 
$k$ is the number of cells in the partition; 
$\cX = [0,1]^p$ with $p=10$. 
There is no randomization of the labels.  
\label{tab.hoRF.dim10}}
\begin{center}
\begin{tabular}{ccc}
\hline\noalign{\smallskip}
& Single tree & Large forest 
\\
\noalign{\smallskip}\hline\noalign{\smallskip}
%
\begin{minipage}{2cm}\centerline{No bootstrap} 
\centerline{$\mathtt{mtry}=p$} \end{minipage} 
& $\displaystyle \frac{0.11}{k^{0.12}} + \frac{1.03 \sigma^{2\vphantom{^2}} k}{n_2}$
& $\displaystyle \frac{0.11}{k^{0.12}} + \frac{1.03 \sigma^{2\vphantom{^2}} k}{n_2}$
\\
\noalign{\smallskip}\hline\noalign{\smallskip}
\begin{minipage}{2cm} \centerline{Bootstrap}  
\centerline{$\mathtt{mtry}=p$} \end{minipage}
& $\displaystyle \frac{0.11}{k^{0.11}} + \frac{1.05 \sigma^{2\vphantom{^2}} k}{n_2}$
& $\displaystyle \frac{0.10}{k^{0.19}} + \frac{0.04 \sigma^{2\vphantom{^2}} k}{n_2}$
\\
\noalign{\smallskip}\hline
\begin{minipage}{2cm} \centerline{No bootstrap}  
\centerline{$\mathtt{mtry}=\lfloor p/3 \rfloor$} \end{minipage}
& $\displaystyle \frac{0.21}{k^{0.18}} + \frac{1.08 \sigma^{2\vphantom{^2}} k}{n_2}$
& $\displaystyle \frac{0.08}{k^{0.25}} + \frac{0.04 \sigma^{2\vphantom{^2}} k}{n_2}$
\\
\noalign{\smallskip}\hline
\begin{minipage}{2cm} \centerline{Bootstrap}  
\centerline{$\mathtt{mtry}=\lfloor p/3 \rfloor$} 
\end{minipage}
& $\displaystyle \frac{0.20}{k^{0.16}} + \frac{1.05 \sigma^{2\vphantom{^2}} k}{n_2}$
& $\displaystyle \frac{0.07}{k^{0.26}} + \frac{0.03 \sigma^{2\vphantom{^2}} k}{n_2}$
\\
\noalign{\smallskip}\hline
\end{tabular} 
\end{center}
\end{table}

\begin{figure}[!ht]
  \begin{center}
  \begin{minipage}{0.45\textwidth}
    \includegraphics[width=\textwidth]{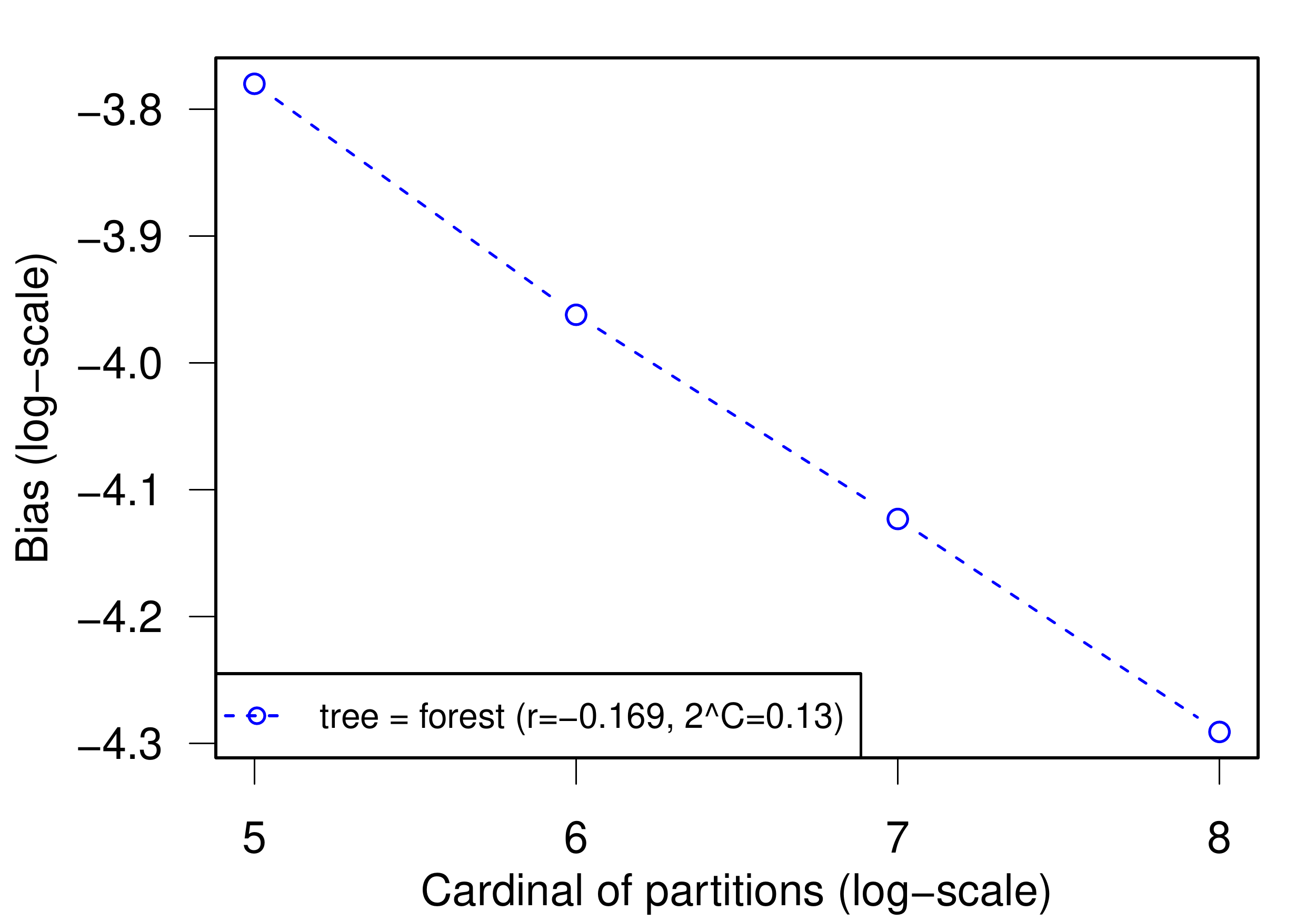}\\
    \centerline{(a) no bootstrap, $\mathtt{mtry} = 5$, $p=5$}
    \end{minipage}
    \hspace{0.01\textwidth}
  \begin{minipage}{0.45\textwidth}
  \includegraphics[width=\textwidth]{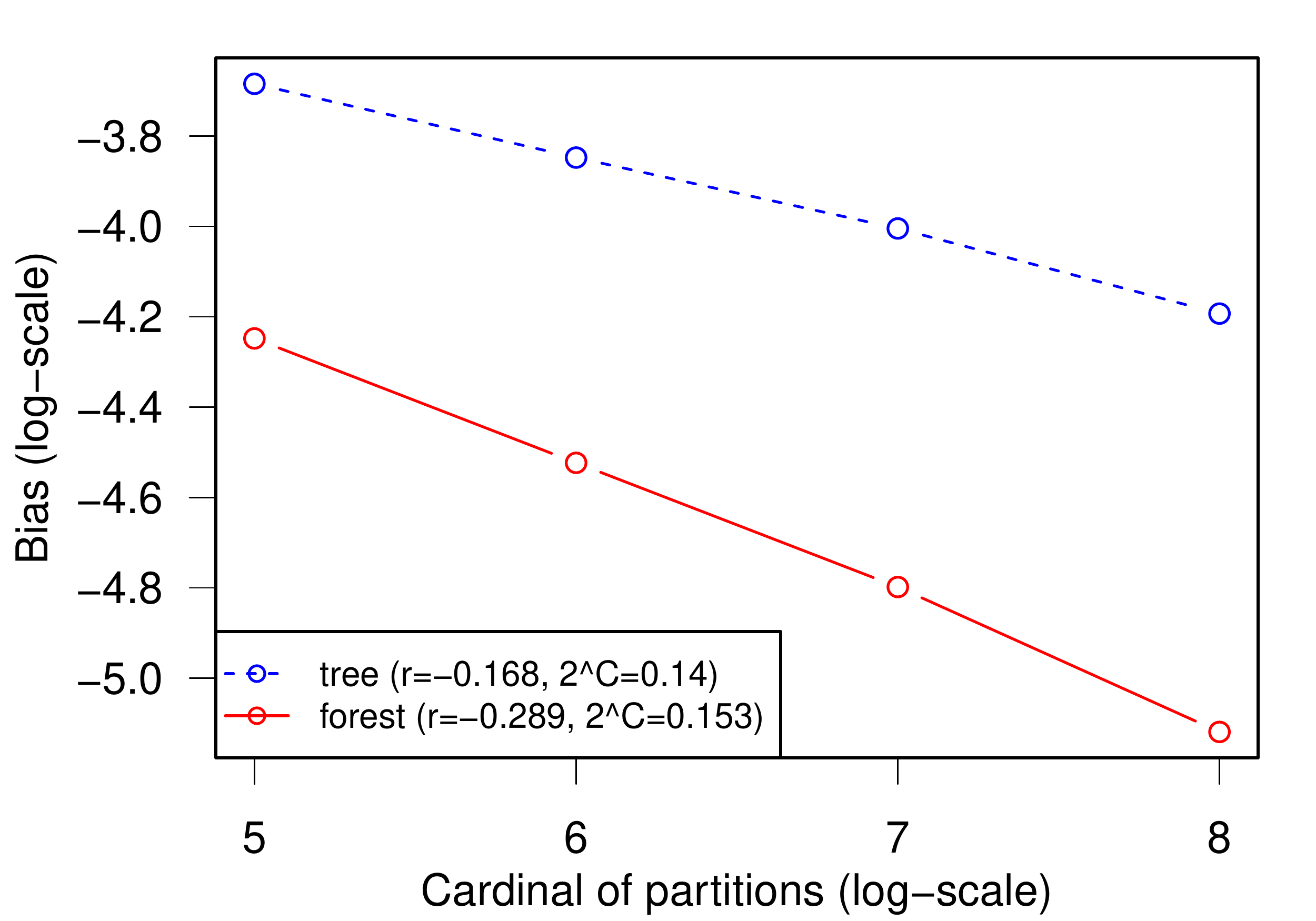}\\
    \centerline{(b) bootstrap, $\mathtt{mtry} = 5$, $p=5$}
  \end{minipage} \\
  \begin{minipage}{0.45\textwidth}
    \includegraphics[width=\textwidth]{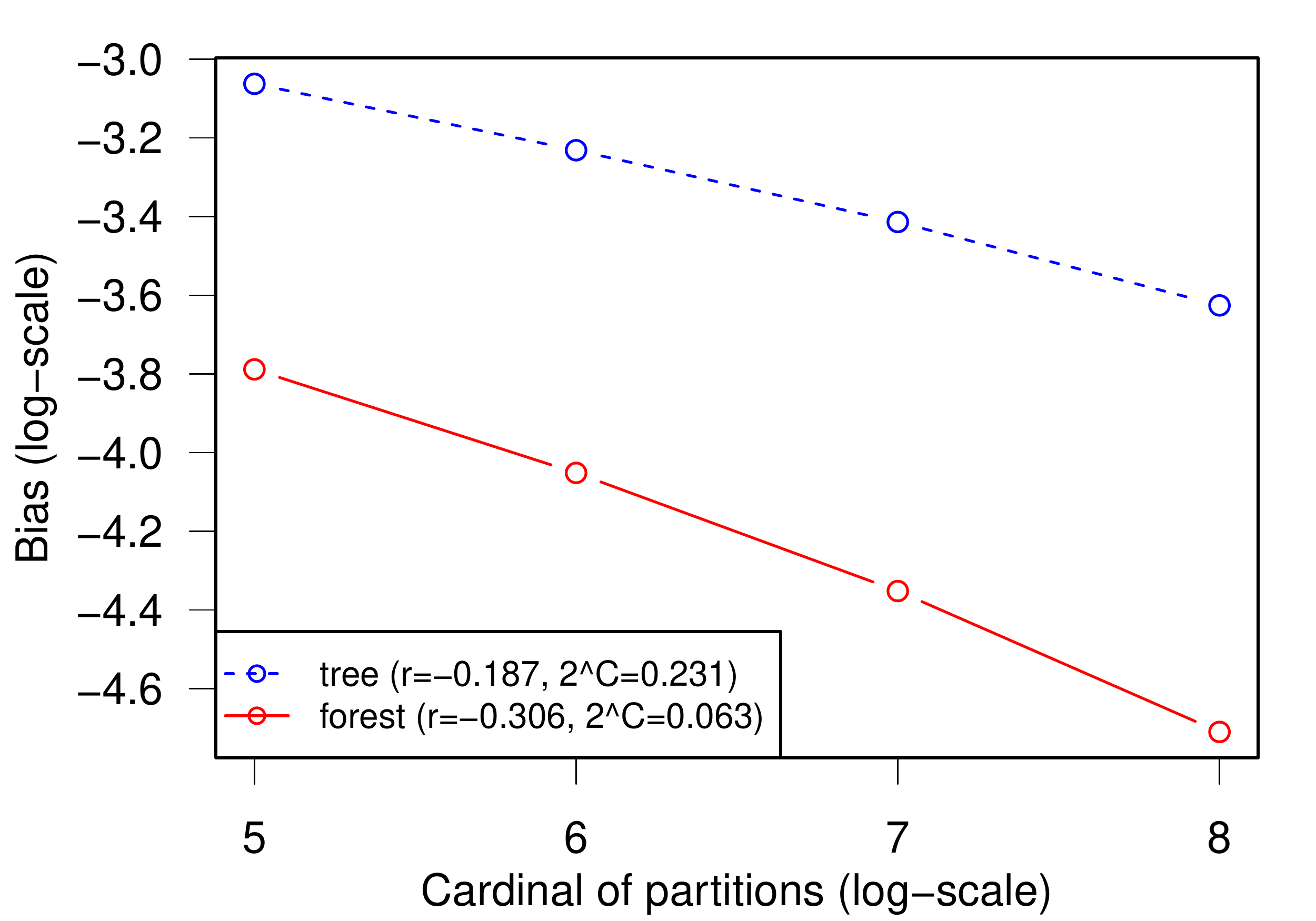}\\
    \centerline{(c) no bootstrap, $\mathtt{mtry} = 1$, $p=5$}
    \end{minipage}
    \hspace{0.01\textwidth}
  \begin{minipage}{0.45\textwidth}
    \includegraphics[width=\textwidth]{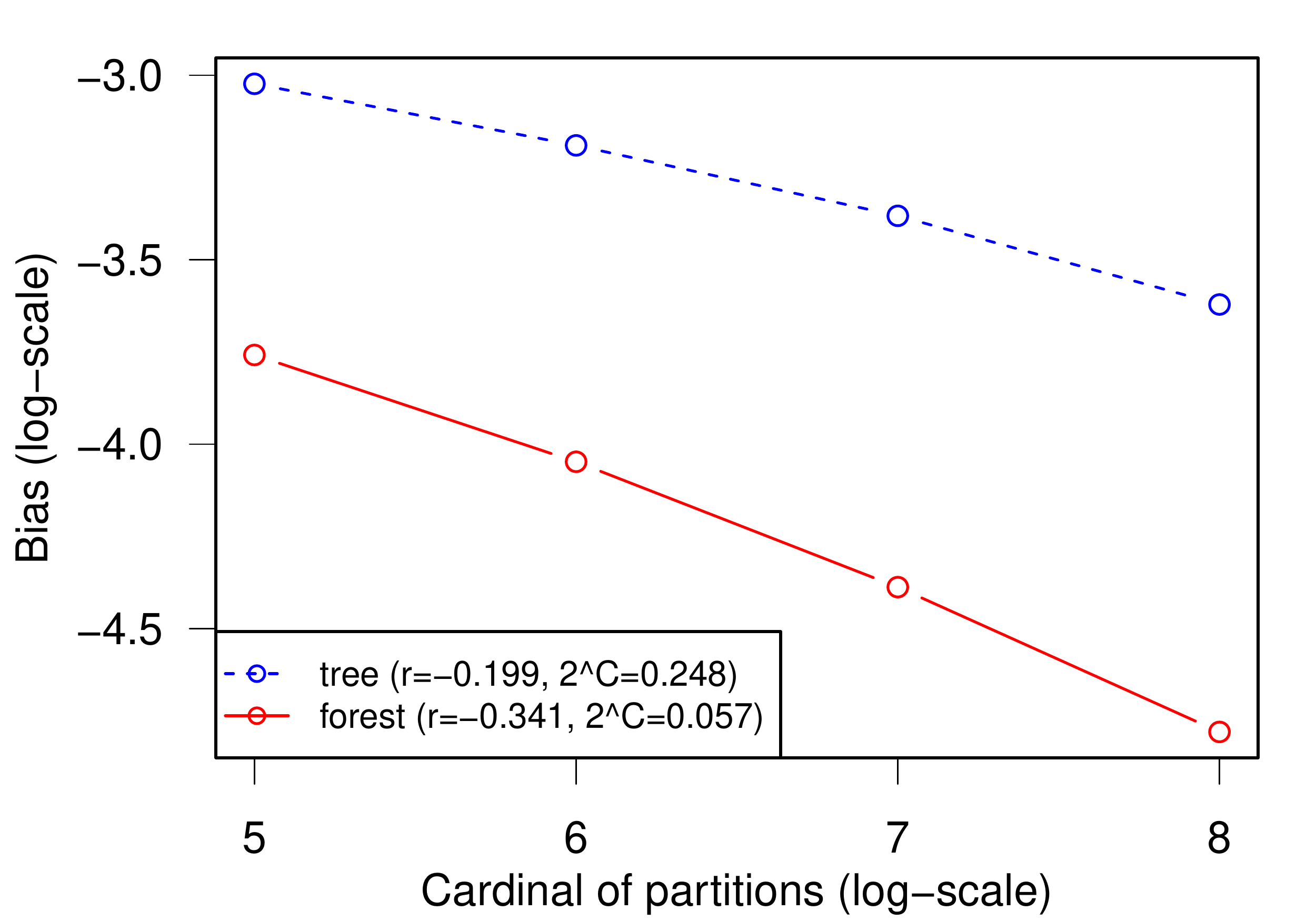}\\
    \centerline{(d) bootstrap, $\mathtt{mtry} = 1$, $p=5$}
    \end{minipage}
    \caption{\label{fig.bias.dim5} 
Estimated values of the approximation error of hold-out trees and 
``large'' forests (in $\log_2$-scale) as a function of the number of leaves (in $\log_2$-scale), 
for the \textbf{Friedman 1} regression function in dimension 
$p=5$, with various values of the parameters 
(bootstrap or not, $\mathtt{mtry} \in \{ p, \lfloor p/3 \rfloor\} $). 
The coefficients $r$ and $C$ respectively denote 
the slope and the intercept of a linear model fitted to the scatter plot. 
}
  \end{center}
\end{figure}

\begin{figure}[!ht]
  \begin{center}
  \begin{minipage}{0.45\textwidth}
    \includegraphics[width=\textwidth]{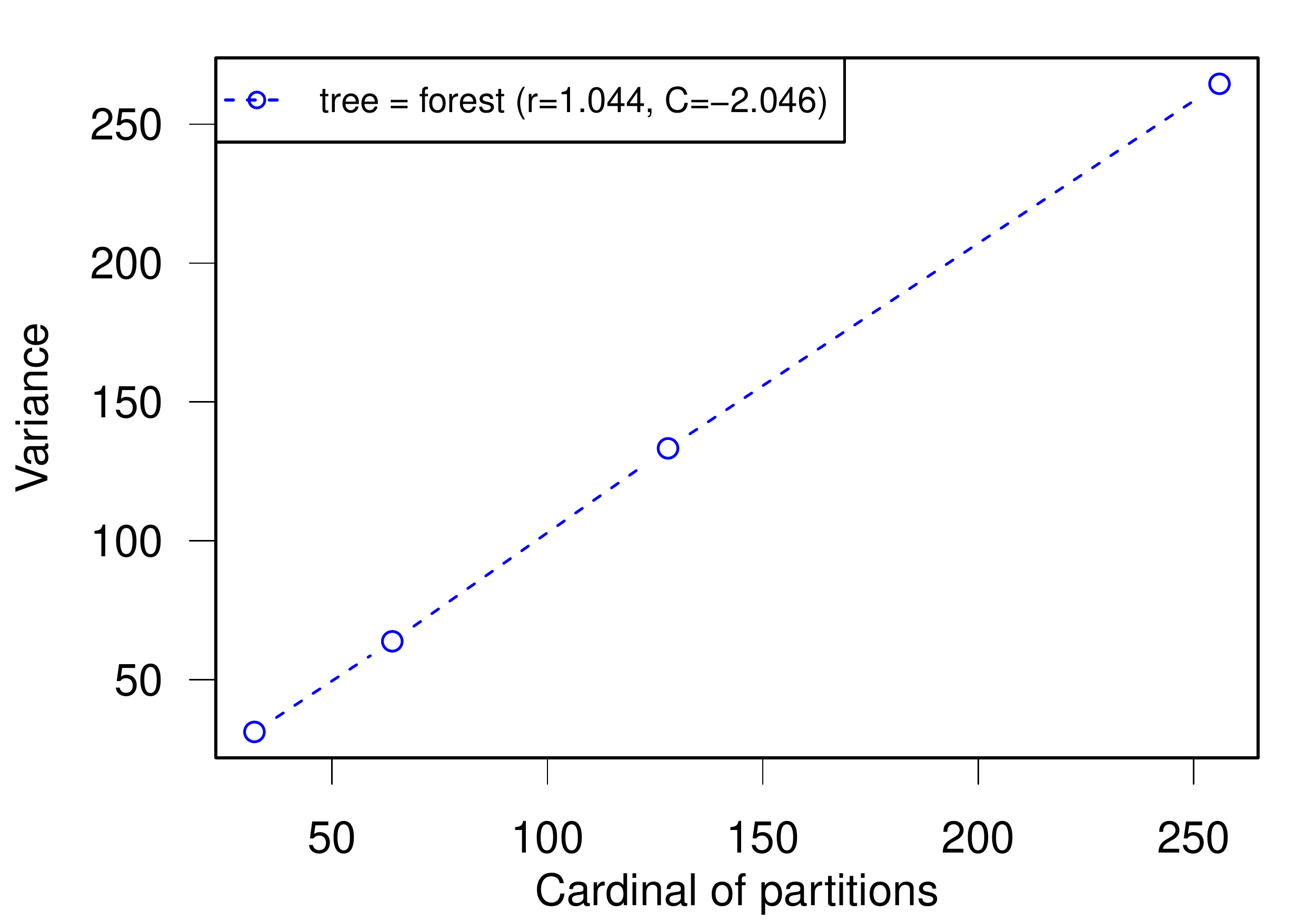}\\
    \centerline{(a) no bootstrap, $\mathtt{mtry} = 5$, $p=5$}
    \end{minipage}
    \hspace{0.01\textwidth}
  \begin{minipage}{0.45\textwidth}
  \includegraphics[width=\textwidth]{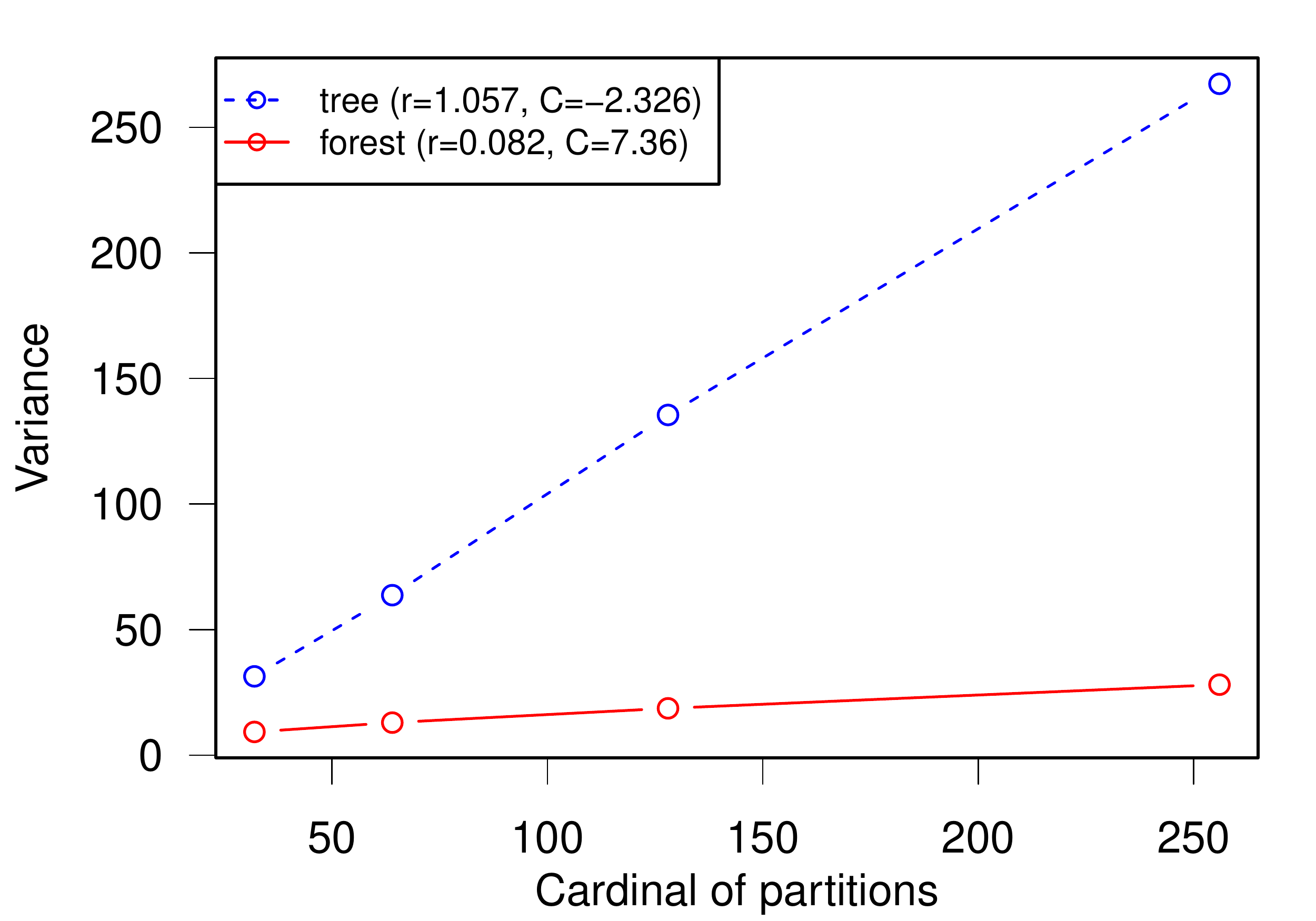}\\
    \centerline{(b) bootstrap, $\mathtt{mtry} = 5$, $p=5$}
  \end{minipage} \\
  \begin{minipage}{0.45\textwidth}
    \includegraphics[width=\textwidth]{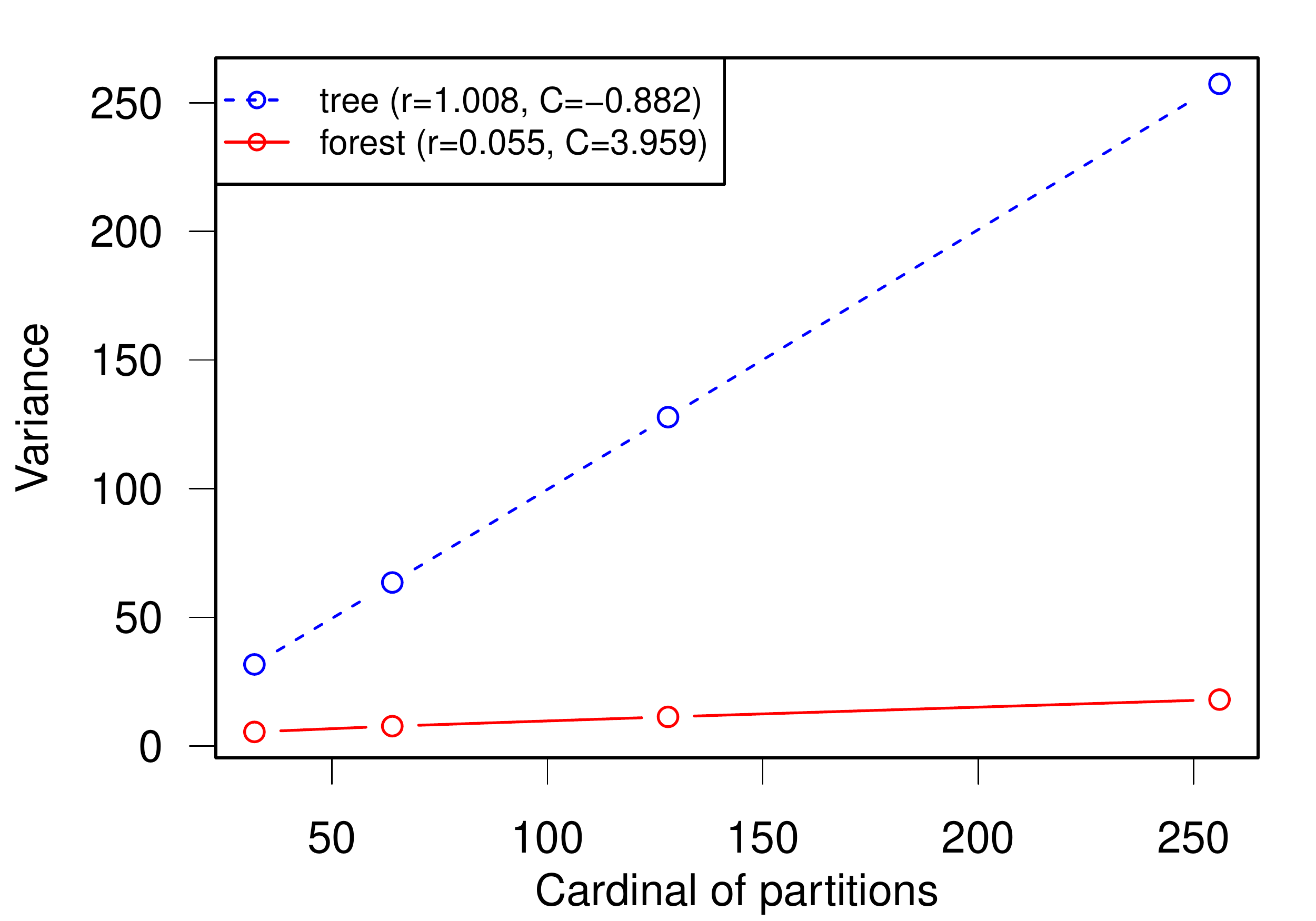}\\
    \centerline{(c) no bootstrap, $\mathtt{mtry} = 1$, $p=5$}
    \end{minipage}
    \hspace{0.01\textwidth}
  \begin{minipage}{0.45\textwidth}
    \includegraphics[width=\textwidth]{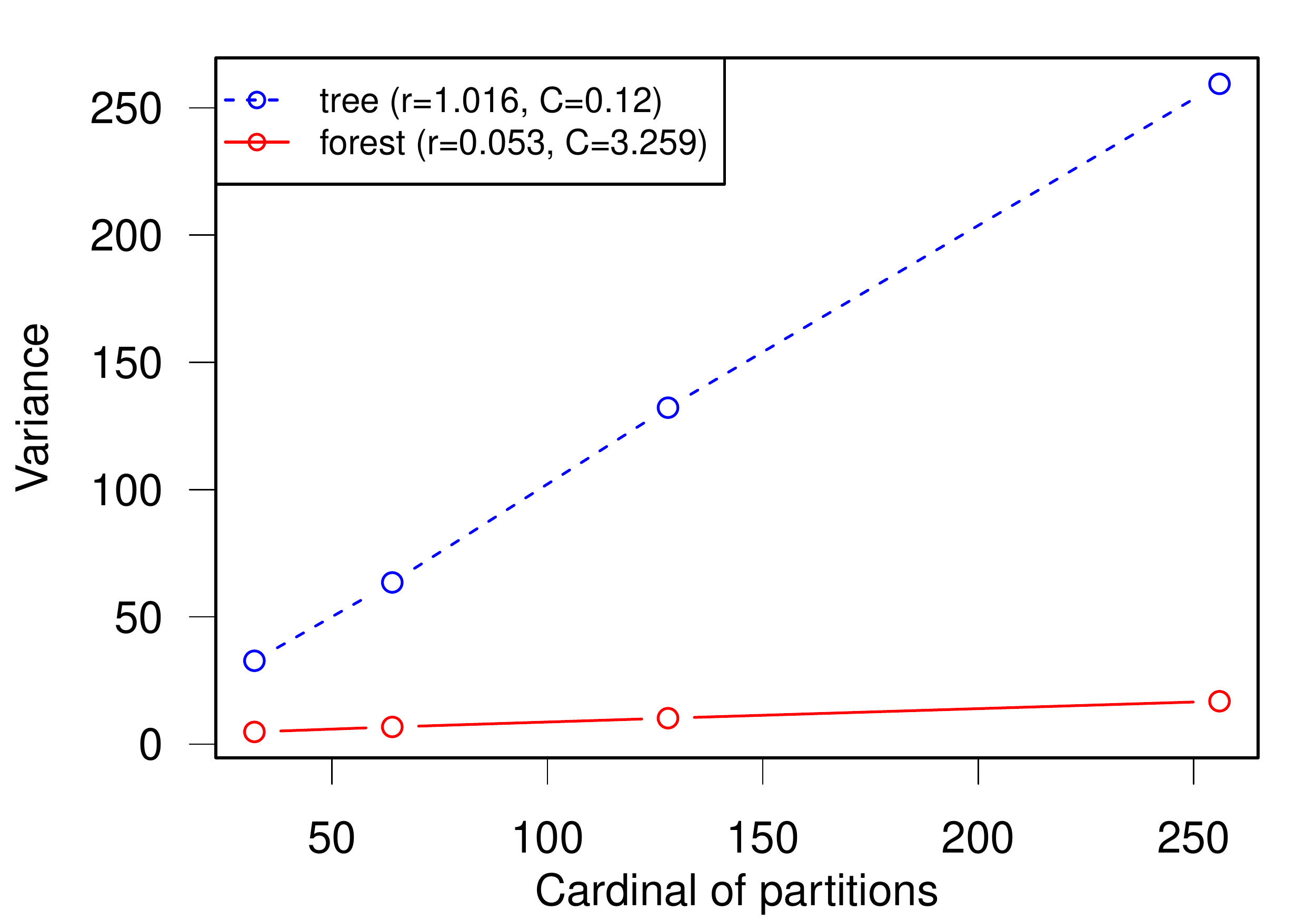}\\
    \centerline{(d) bootstrap, $\mathtt{mtry} = 1$, $p=5$}
    \end{minipage}
    \caption{\label{fig.var.dim5} 
Estimated values of the estimation error (multiplied by $n_2/\sigma^2$) of hold-out trees and 
``large'' forests as a function of the number of leaves, 
for the \textbf{Friedman 1} regression function in dimension 
$p=5$, with various values of the parameters 
(bootstrap or not, $\mathtt{mtry} \in \{ p, \lfloor p/3 \rfloor\} $). 
The coefficients $r$ and $C$ respectively denote 
the slope and the intercept of a linear model fitted to the scatter plot. 
}
  \end{center}
\end{figure}

\begin{figure}[!ht]
  \begin{center}
  \begin{minipage}{0.45\textwidth}
    \includegraphics[width=\textwidth]{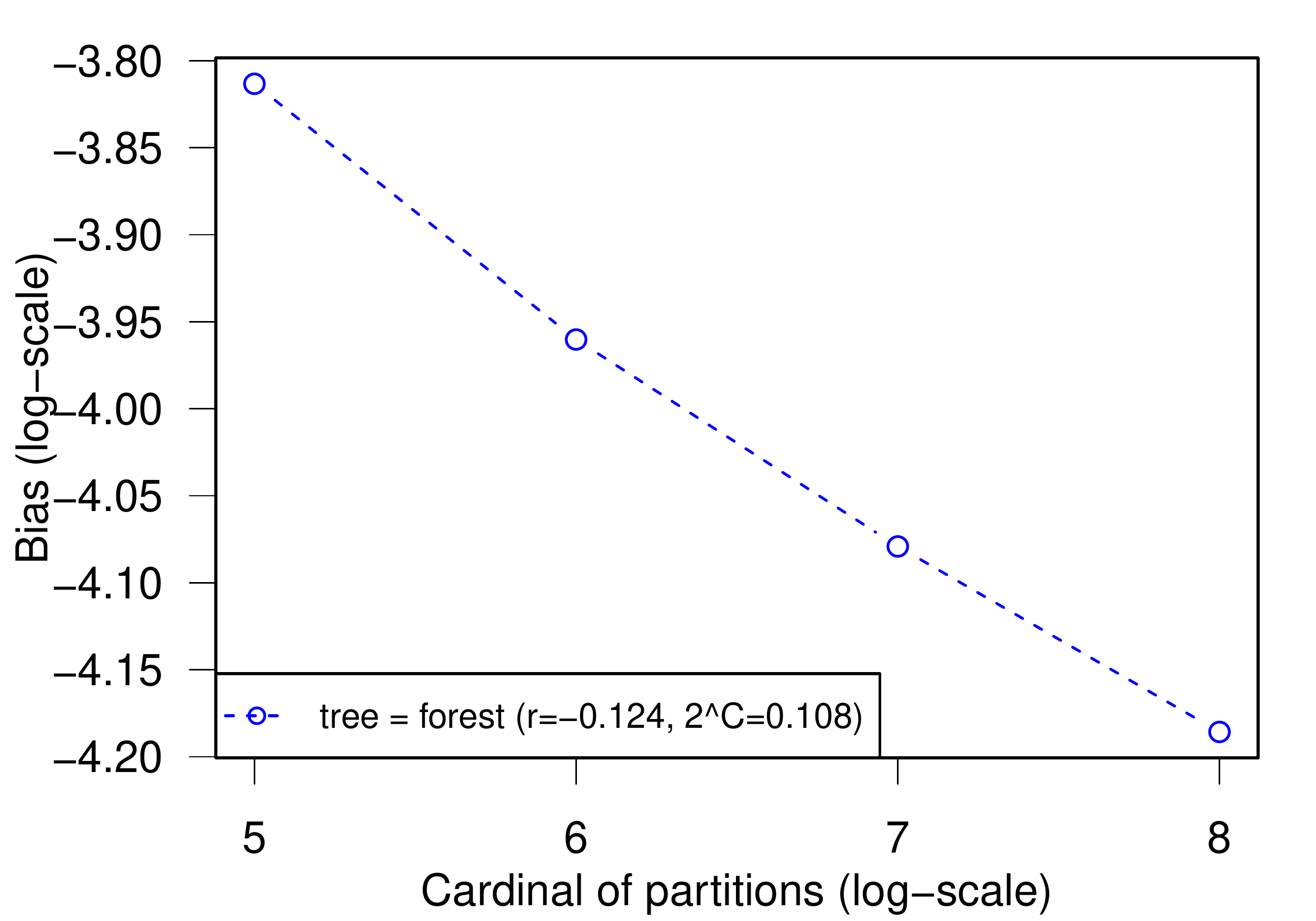}\\
    \centerline{(a) no bootstrap, $\mathtt{mtry} = 10$, $p=10$}
    \end{minipage}
    \hspace{0.01\textwidth}
  \begin{minipage}{0.45\textwidth}
  \includegraphics[width=\textwidth]{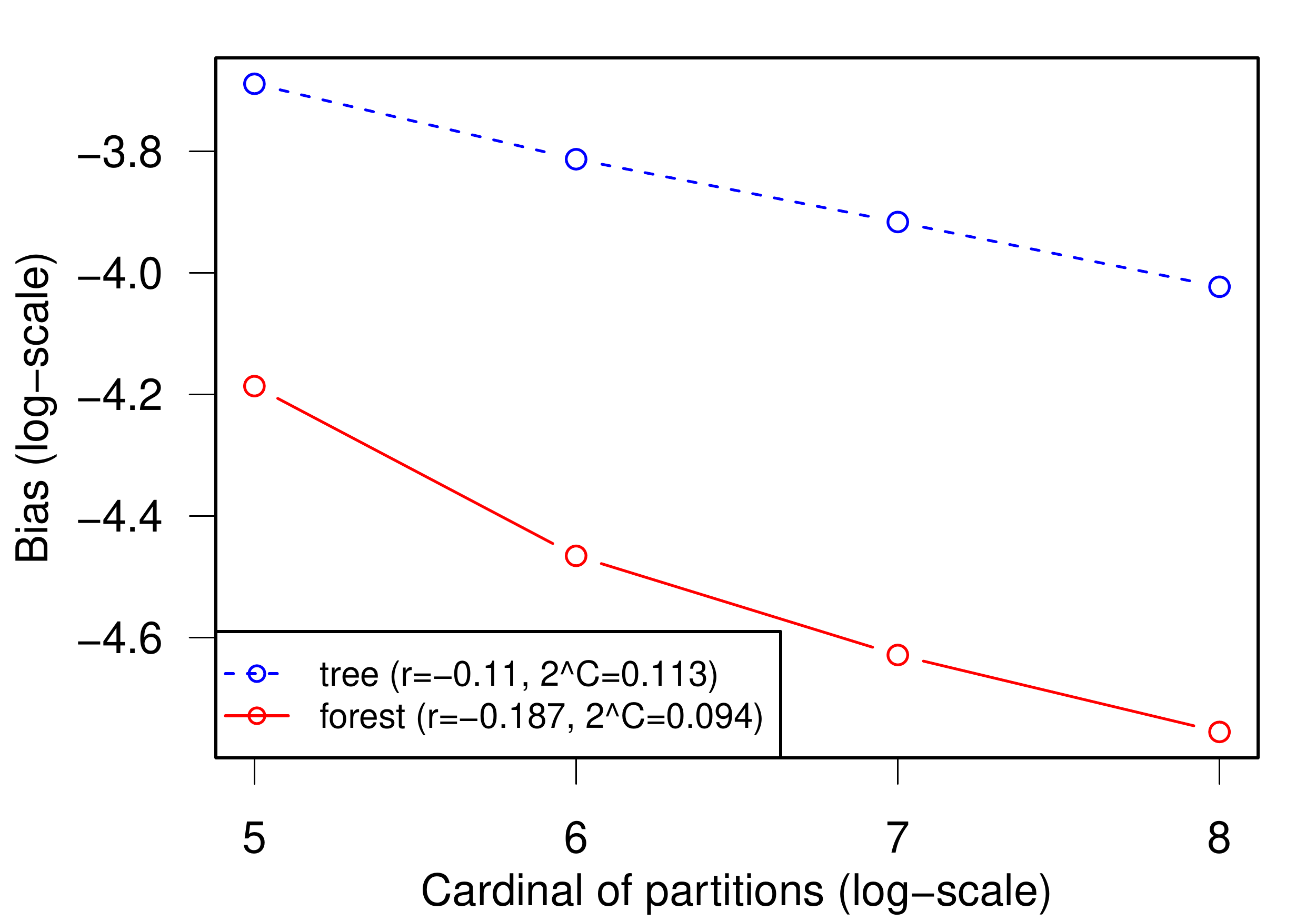}\\
    \centerline{(b) bootstrap, $\mathtt{mtry} = 10$, $p=10$}
  \end{minipage} \\
  \begin{minipage}{0.45\textwidth}
    \includegraphics[width=\textwidth]{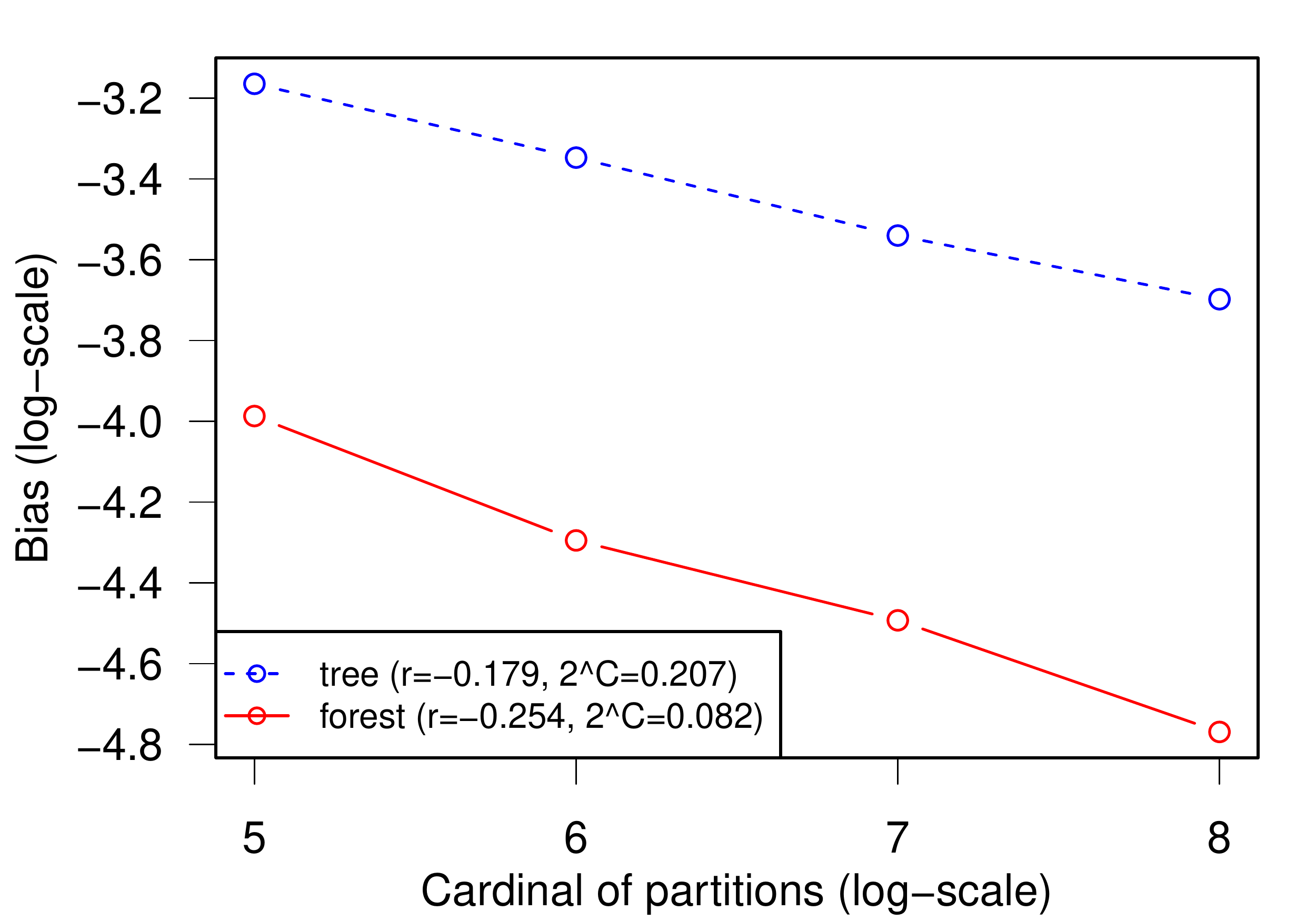}\\
    \centerline{(c) no bootstrap, $\mathtt{mtry} = 3$, $p=10$}
    \end{minipage}
    \hspace{0.01\textwidth}
  \begin{minipage}{0.45\textwidth}
    \includegraphics[width=\textwidth]{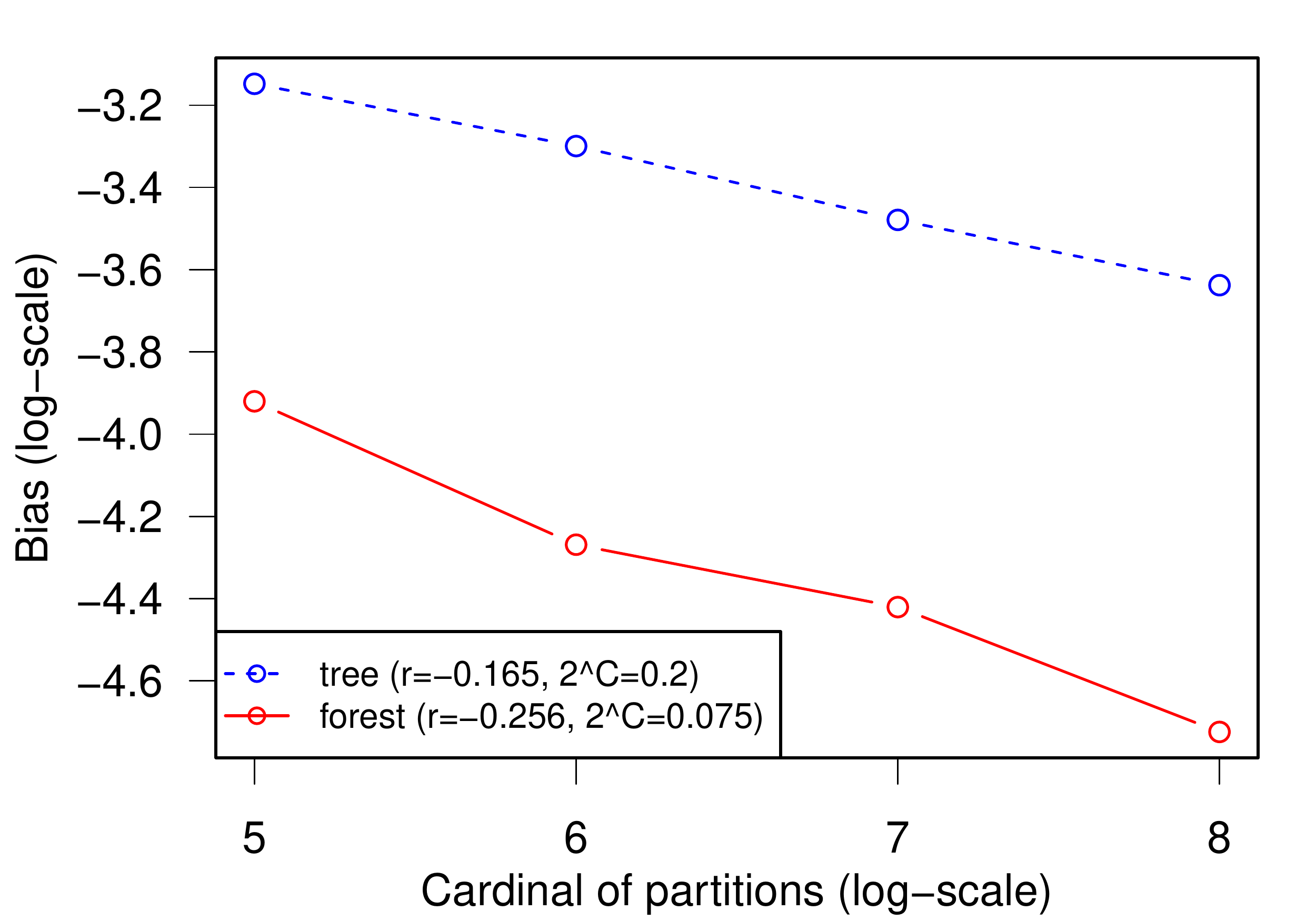}\\
    \centerline{(d) bootstrap, $\mathtt{mtry} = 3$, $p=10$}
    \end{minipage}
    \caption{\label{fig.bias.dim10} 
Estimated values of the approximation error of hold-out trees and 
``large'' forests (in $\log_2$-scale) as a function of the number of leaves (in $\log_2$-scale), 
for the \textbf{Friedman 1} regression function in dimension 
$p=10$, with various values of the parameters 
(bootstrap or not, $\mathtt{mtry} \in \{ p, \lfloor p/3 \rfloor\} $). 
The coefficients $r$ and $C$ respectively denote 
the slope and the intercept of a linear model fitted to the scatter plot. }
  \end{center}
\end{figure}

\begin{figure}[!ht]
  \begin{center}
  \begin{minipage}{0.45\textwidth}
    \includegraphics[width=\textwidth]{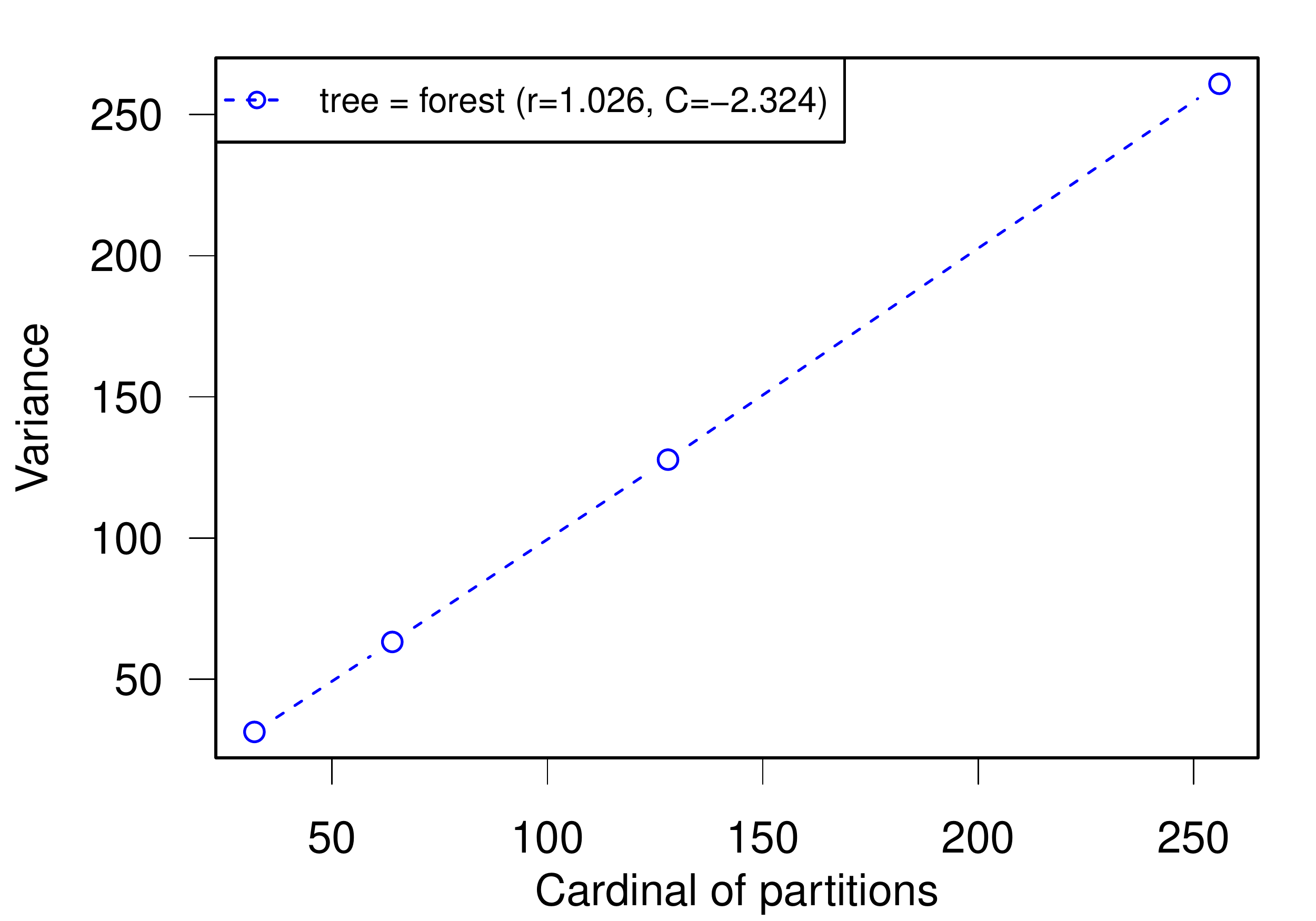}\\
    \centerline{(a) no bootstrap, $\mathtt{mtry} = 10$, $p=10$}
    \end{minipage}
    \hspace{0.01\textwidth}
  \begin{minipage}{0.45\textwidth}
  \includegraphics[width=\textwidth]{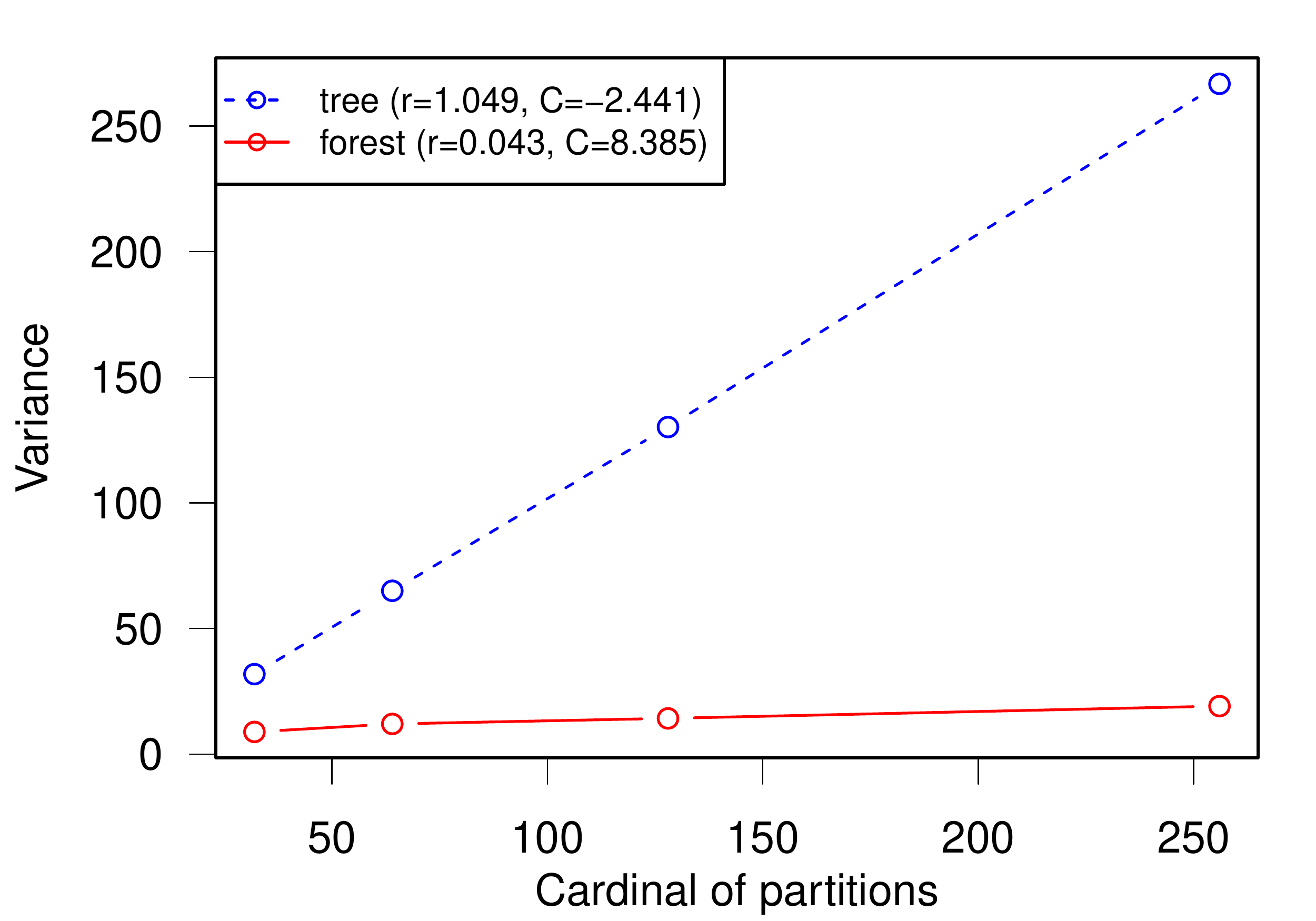}\\
    \centerline{(b) bootstrap, $\mathtt{mtry} = 10$, $p=10$}
  \end{minipage} \\
  \begin{minipage}{0.45\textwidth}
    \includegraphics[width=\textwidth]{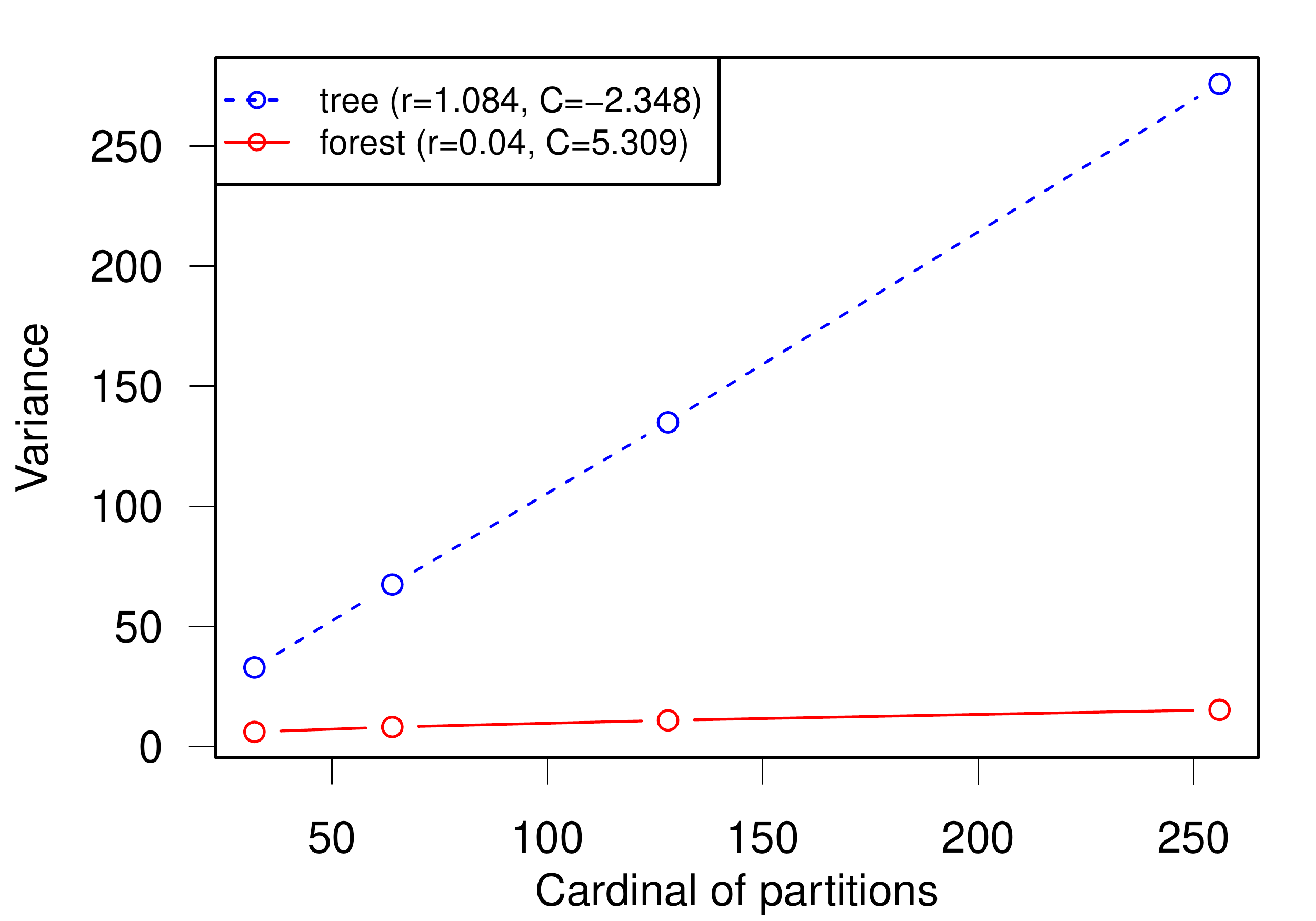}\\
    \centerline{(c) no bootstrap, $\mathtt{mtry} = 3$, $p=10$}
    \end{minipage}
    \hspace{0.01\textwidth}
  \begin{minipage}{0.45\textwidth}
    \includegraphics[width=\textwidth]{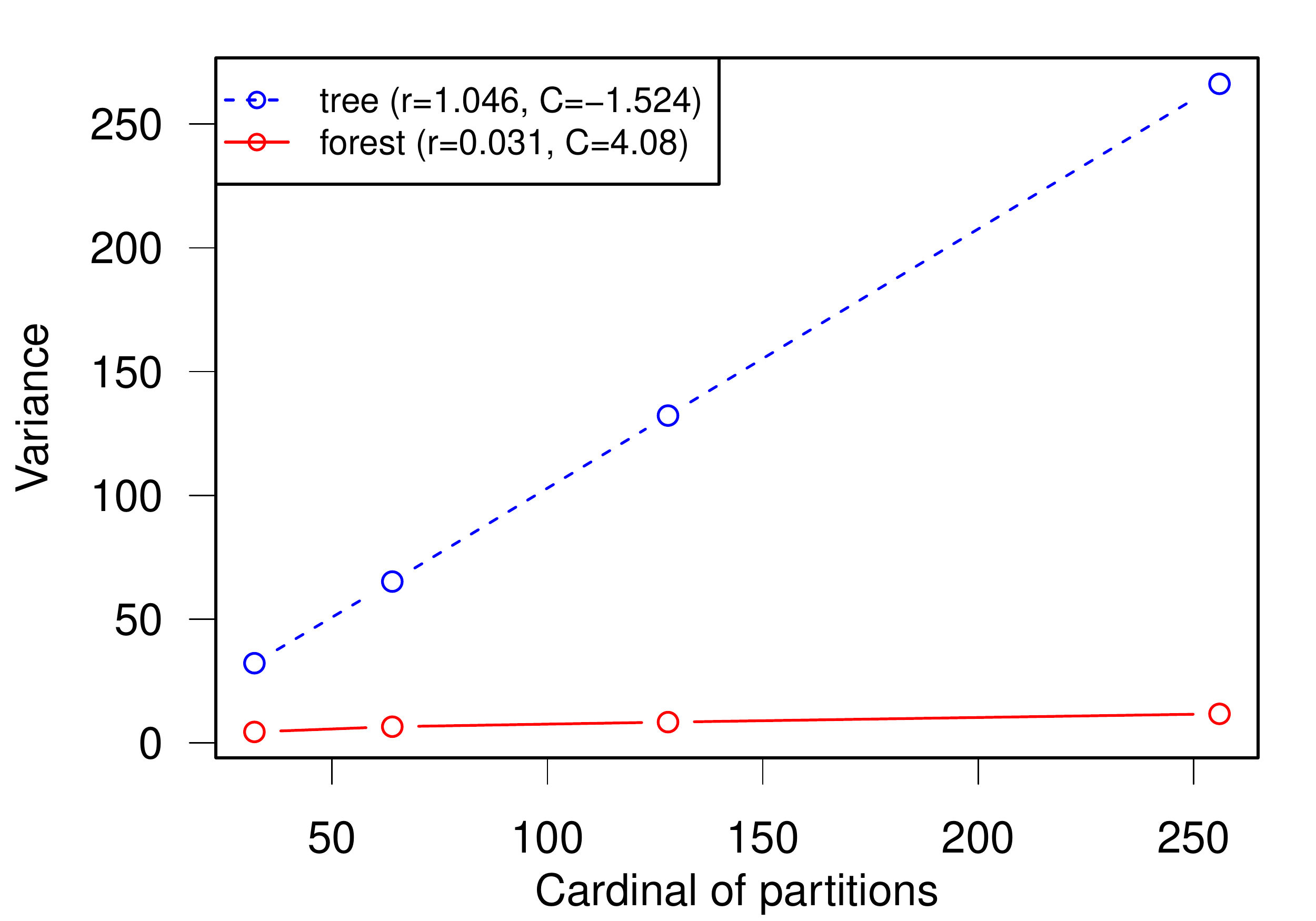}\\
    \centerline{(d) bootstrap, $\mathtt{mtry} = 3$, $p=10$}
    \end{minipage}
    \caption{\label{fig.var.dim10} 
Estimated values of the estimation error (multiplied by $n_2/\sigma^2$) of hold-out trees and 
``large'' forests as a function of the number of leaves, 
for the \textbf{Friedman 1} regression function in dimension 
$p=10$, with various values of the parameters 
(bootstrap or not, $\mathtt{mtry} \in \{ p, \lfloor p/3 \rfloor\} $). 
The coefficients $r$ and $C$ respectively denote 
the slope and the intercept of a linear model fitted to the scatter plot. 
}
  \end{center}
\end{figure}

\end{document}